\documentclass[11pt]{article} 

\usepackage{tikz}
\usetikzlibrary{positioning}
\usepackage[utf8]{inputenc} 
\usepackage{amsmath}
\usepackage{amsmath,breqn}
\usepackage{amssymb}
\usepackage{graphicx}
\usepackage{epsfig,epstopdf}
\usepackage{pdflscape}
\usepackage{rotating,lscape,afterpage}
\usepackage{amsthm}

\usepackage{geometry} 
\usepackage{graphicx} 

\usepackage{booktabs} 
\usepackage{booktabs} 
\usepackage{array} 
\usepackage{paralist} 
\usepackage{verbatim} 
\usepackage{subfig} 

\usepackage{fancyhdr} 
\usepackage{tikz}
\usetikzlibrary{positioning}
\usepackage[utf8]{inputenc} 
\usepackage{amsmath}
\usepackage{amsmath,breqn}
\usepackage{amssymb}
\usepackage{graphicx}
\usepackage{epsfig,epstopdf}
\usepackage{pdflscape}
\usepackage{rotating,lscape,afterpage}
\usepackage{amsthm}

\usepackage{geometry} 
\usepackage{graphicx} 

\usepackage{sectsty}
\allsectionsfont{\sffamily\mdseries\upshape} 

\usepackage[nottoc,notlof,notlot]{tocbibind} 
\usepackage[titles,subfigure]{tocloft} 

\numberwithin{equation}{section}

\begin{document}

\title{The evolution problem for the 1D nonlocal Fisher-KPP equation with a top hat kernel. Part 2. The Cauchy problem on a finite interval}

\author{D. J. Needham\\School of Mathematics\\ University of Birmingham\\B15 2TT UK\\d.j.needham@bham.ac.uk\and J. Billingham\\School of Mathematical Sciences\\ University of Nottingham\\
NG7 2RD UK\\john.billingham@nottingham.ac.uk}

\maketitle
\vspace{10pt}

\begin{abstract}
In the second part of this series of papers, we address the same Cauchy problem that was considered in part 1 (see \cite{NBLM}), namely the nonlocal Fisher-KPP equation in one spatial dimension,
\[
u_t = D u_{xx} + u(1-\phi*u),
\]
where $\phi*u$ is a spatial convolution with the top hat kernel, $\phi(y) \equiv H\left(\frac{1}{4}-y^2\right)$, except that now the spatial domain is the finite interval $[0,a]$ rather than the whole real line. Consequently boundary conditions are required at the interval end-points, and we address the situations when these boundary conditions are of either Dirichlet or Neumann type. This model forms a natural extension to the classical Fisher-KPP model, with the introduction of the simplest possible nonlocal effect into the saturation term. Nonlocal reaction-diffusion models arise naturally in a variety of (frequently biological or ecological) contexts, and as such it is of fundamental interest to examine its properties in detail, and to compare and contrast these with the well known properties of the classical Fisher-KPP model.
\end{abstract}

\section{Introduction}\label{sec_intro}
In this second paper we consider the evolution problem detailed in part 1 of this pair of papers \cite{NBLM} (and henceforth referred to as (NB1)) for the nonlocal Fisher-KPP equation with top hat kernel, but now on a finite spatial interval $[0,a]$, where $a$ is the dimensionless interval length, scaled relative to the nonlocal length scale. As described in (NB1) and references therein, nonlocal reaction-diffusion equations arise in many different scientific areas, and finite domain effects can be relevant in many of them, particularly in biomedical applications, \cite{VolpBio}. The reader is referred to the introduction of (NB1) for a more general discussion, which need not be repeated here. However, it is worth emphasising that the key feature of the inclusion of the nonlocal term is the introduction of a nonlocal length scale into the model. There are thus two dimensionless parameters in the model, namely $D$, which measures the square of the ratio of the diffusion length scale (based on the kinetic time scale) to the nonlocal length scale, and $a$, which measures the ratio of the domain length scale to the nonlocal length scale. We examine both the Dirichlet and the Neumann models for $(a,D)$ throughout the positive quadrant of the parameter plane. As may be anticipated, significant structural differences in behaviour between the classical Fisher-KPP model and this natural nonlocal extended model become evermore present as the parameter $D$ decreases, and these principal, and marked, differences are fully explored and exposed. 

The simplicity of the top hat kernel allows us to develop both a qualitative and quantitative, asymptotic description of the structure and dynamics of the solutions, whose qualitative features (hump formation and multiple stable steady state attractors), which have previously been unidentified, are expected to remain relevant for other localised kernels that are less amenable to this type of analysis. In particular, we study regimes far away from linear and weakly nonlinear limits, and address the fully nonlinear interaction between nonlocality, diffusion and competition in the underlying nonlocal PDE.

We study the evolution problem, firstly with Dirichlet and secondly with Neumann boundary conditions at the two ends of the interval. The evolution problem is as follows,
\begin{equation}
    u_t = Du_{xx} + u\left(1 - \int_{\alpha(x)}^{\beta(x)}{u(y,t)dy}\right)~~\text{with}~~(x,t)\in D_{\infty}, \label{eqn1.1}
\end{equation}
where, with $0<a\le1/2,$ we have,
\begin{equation}
  \alpha(x) = 0~~\text{and}~~\beta(x)=a~~\forall~~x\in[0,a], \label{eqn1.2}  
\end{equation}
whilst, with $a>1/2,$ we have,
\begin{equation}
    \alpha(x) = 
    \begin{cases}
    0,~0\le x\le \frac{1}{2},\\
    x-\frac{1}{2},~\frac{1}{2}<x\le a,
    \end{cases}
    \label{eqn1.3}
\end{equation}
and,
\begin{equation}
    \beta(x) = 
    \begin{cases}
    x+\frac{1}{2},~0\le x <a-\frac{1}{2},\\
    a,~a-\frac{1}{2}\le x \le a.
    \end{cases}
    \label{eqn1.4}
\end{equation}
The associated initial condition is,
\begin{equation}
    u(x,0)=u_0(x),~x\in[0,a], \label{eqn1.5}
\end{equation}
whilst the boundary conditions are, for the Dirichlet problem,
\begin{equation}
    u(0,t) = u(a,t) = 0,~t\in(0,\infty), \label{eqn1.6}
\end{equation}
or for the Neumann problem,
\begin{equation}
    u_x(0,t)=u_x(a,t)=0,~t\in(0,\infty).
    \label{eqn1.7}
\end{equation}
In the above, $D$ is the constant diffusion coefficient, and the initial data has $u_0\in C([0,a])\cap PC^1([0,a])$, is non-negative and nontrivial. Also,
\begin{equation}
    D_{\infty}=(0,a)\times(0,\infty). \label{eqn1.8}
\end{equation}
We henceforth refer to the Dirichlet problem as (DIVP) and the Neumann problem as (NIVP), and throughout we will consider classical solutions to these two evolution problems. 

To begin with, it is a straightforward consequence of the parabolic strong maximum principle (and non-negative, nontrivial initial data) that, with $u:D_{\infty}\to\mathbb{R}$ being a solution to either of (DIVP) or (NIVP), then,
\begin{equation}
    u(x,t) > 0~~\forall~~(x,t)\in D_{\infty}. \label{eqn1.9}
\end{equation}
It is then a consequence of the parabolic comparison theorem (with parabolic operator $N[u] \equiv u_t - Du_{xx} - u)$ that,
\begin{equation}
    u(x,t) \le Me^t~~\forall~~(x,t)\in \overline{D} _{\infty}, \label{eqn1.10}
\end{equation}
with $M=$ sup$_{y\in[0,a]}u_0(y)$. The a priori bounds in (\ref{eqn1.9}) and (\ref{eqn1.10}) readily guarantee the existence and uniqueness of a classical solution to either of (DIVP) and (NIVP) (without repetition, this follows very closely from the approach laid out in section 2 of (NB1) for the Cauchy problem on the real line). A further simple application of the parabolic comparison theorem establishes that, for the solution to (DIVP), there exists a positive constant $A$, for which $u_0(x)\le A\sin\left(\frac{\pi x}{a}\right)~~ \forall~~ x\in [0,a]$, and such that,
\begin{equation}
    u(x,t) \le Ae^{\left(1-\frac{{\pi}^2D}{a^2}\right)t}\sin\left(\frac{\pi x}{a}\right)~~\forall~~(x,t)\in \overline{D}_{\infty}. \label{eqn1.11}
\end{equation}
  The remainder of the paper is structured as follows. Section 2 addresses (DIVP) in detail. Exact solutions are constructed for $0<a\le 1/2$ which allow for a complete analysis of (DIVP). For $a>1/2$ we first consider the existence of nontrivial and non-negative steady states associated with (DIVP). We achieve this by fixing $D$ as an unfolding parameter and regarding $a$ as a bifurcation parameter. The detailed bifurcation structure is examined numerically via pseudo-arclength continuation (see Appendix~\ref{app_num} or \cite{arclength}), with the asymptotic structure of steady states examined in detail in a number of relevant limits in the positive quadrant of the $(a,D)$ parameter plane. The temporal linear stability of these steady states is examined by analysis of the associated linear self-adjoint nonlocal eigenvalue problem. This then enables us to make conjectures relating to the large-$t$ attractors for (DIVP), and these conjectures are then supported by careful numerical solutions to (DIVP). In section 3 we follow a similar approach to developing a detailed analysis to (NIVP), and compare and contrast the results with the corresponding results in section 2. We draw some relevant conclusions in section 4. Descriptions of the numerical methods used in the paper are given in Appendix~\ref{app_num}.
  
  \section{The Dirichlet problem (DIVP)}
 In this section we focus on the Dirichlet problem (DIVP). To begin with, we observe that when the initial data is trivial, the solution to the associated (DIVP) is the equilibrium solution
 \begin{equation}
     u(x,t)=u_e=0~~\forall~~(x,t)\in \bar{D}_{\infty}.  \label{eqn2.1}
 \end{equation} 
 For initial data with $\parallel$$u_0$$\parallel_{\infty}$ small, it is straightforward to develop a linearised theory for (DIVP). For brevity we do not present the details, but observe that this linearised theory for (DIVP) establishes that the trivial equilibrium soluton is locally temporally asymptotically stable when $D>a^2/\pi^2$ but temporally unstable when $0<D<a^2/\pi^2$. These conclusions can be extended, via (\ref{eqn1.9}) and (\ref{eqn1.11}), to establish that the trivial equilibrium solution is globally temporally asymptotically stable when $D>a^2/\pi^2$, and at least temporally Liapunov stable when $D=a^2/\pi^2.$ Thus, when $D>a^2/\pi^2$, the solution to (DIVP) decays to zero exponentially in $t$ as $t\to\infty$, uniformly for $x\in[0,a]$. The remainder of this section focuses on considering the nature of the solution to (DIVP) when $0<D<a^2/\pi^2$, and in particular the structure of the solution for $t$ large.To continue, it is convenient to consider the cases when $0<a\le 1/2$ and $a>1/2$ separately.
 
 \subsection{$0<a\le 1/2$}
 
\begin{figure}
\begin{center}
\includegraphics[width=0.8\textwidth]{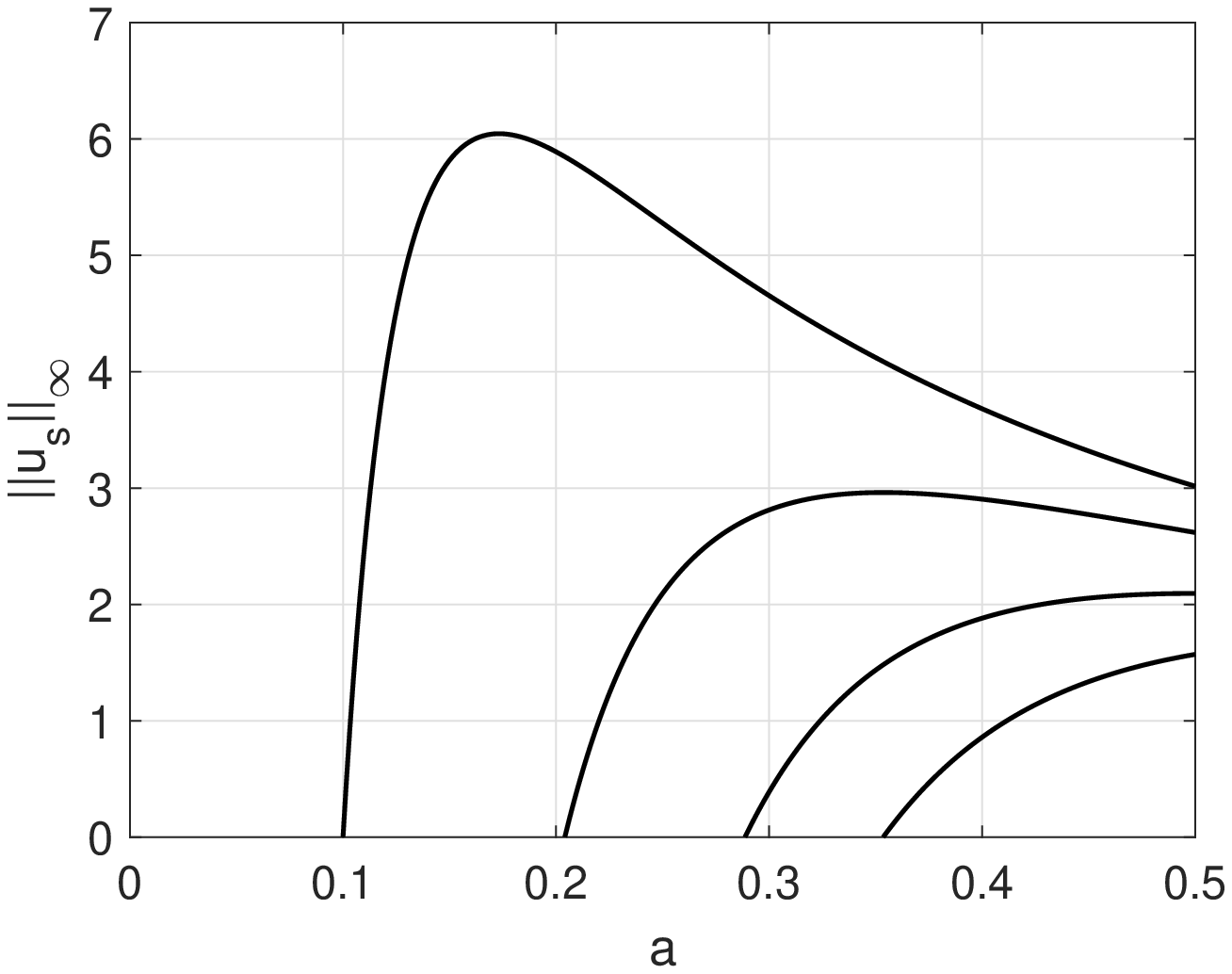}
\caption{$||u_s||_{\infty}$ as a function of $a$ for $\pi^2 D = \frac{1}{100}$, $\frac{1}{24}$, $\frac{1}{12}$ and $\frac{1}{8}$.}\label{fig3}
\end{center}
\end{figure}

 In this case it follows from (\ref{eqn1.2}), that equation (\ref{eqn1.1}) becomes,
 \begin{equation}
     u_t = Du_{xx} + u(1-a\bar{u}(t)),~~(x,t)\in D_{\infty}, \label{eqn2.2}
 \end{equation}
 where $\bar{u}:[0,\infty)\to \mathbb{R}$ is the spatial mean value of $u$ over the interval $[0,a],$ given by,
 \begin{equation}
     \bar{u}(t) = \frac{1}{a} \int_0^a{u(y,t)}dy~~\forall~~t\in[0,\infty). \label{eqn2.3}
 \end{equation}
 Before considering (DIVP) in detail, it is instructive to first consider the possible non-negative and nontrivial steady states which satisfy equation (\ref{eqn2.2}) and boundary conditions (\ref{eqn1.6}), as we anticipate that the solution to (DIVP) will be attracted to one such steady state as $t\to\infty$ for $0<D<a^2/\pi^2$. A non-negative steady state $u_s:[0,a]\to\mathbb{R}$ to (DIVP) satisfies,
 \begin{equation}
     Du_s'' + u_s(1-a\bar{u}_s) = 0, ~~x\in(0,a),
      \label{eqn2.4'}
 \end{equation}
 with,
 \begin{equation}
     u_s(x)\ge0~~\forall~~x\in[0,a], \label{eqn2.5'}
 \end{equation}
 subject to the boundary conditions,
 \begin{equation}
     u_s(0)=u_s(a)=0,  \label{eqn2.6'}
 \end{equation}
 and with the constant $\bar{u}_s$ given by,
 \begin{equation}
     \bar{u}_s = \frac{1}{a}\int_0^a{u_s(y)}dy. \label{2.7'}
 \end{equation}
 We will refer to this boundary value problem as (BVP). It is straightforward to analyse (BVP) directly. Apart from the trivial solution, (BVP) has an  additional nontrivial solution if and only if $0<D<a^2/\pi^2$, and this solution is given by,
 \begin{equation}
     u_s(x) = \frac{\pi}{2a}\left(1 - \frac{\pi^2D}{a^2}\right)\sin\left(\frac{\pi x}{a}\right)~~\forall~~x\in[0,a],  \label{eqn2.8'}
 \end{equation}
 with,
 \begin{equation}
     \bar{u}_s = \frac{1}{a}\left(1 - \frac{\pi^2D}{a^2}\right) = \frac{1}{a}||u_s||_1 = \frac{2}{\pi}||u_s||_{\infty}. \label{eqn2.9'}
 \end{equation}
 We observe that the nontrivial steady state is a positive single hump, symmetric about $x=\frac{1}{2}a$, and spanning the whole interval. It is convenient to consider this steady state with $0<D<1/4\pi^2$ fixed, and allowing $a$ to vary in the interval of existence $(\pi\sqrt{D},1/2].$ From (\ref{eqn2.8'}) and (\ref{eqn2.9'}), we see that $u_s$ emerges at a steady state transcritical bifurcation from the equilibrium state $u_e=0$, when $a=\pi\sqrt{D}$ (noting that the nontrivial steady solution generated by this bifurcation for $a\in(0,\pi\sqrt{D})$ is negative, and thus ruled out) and develops for $a\in(\pi\sqrt{D},1/2].$ For $1/12\pi^2\le D<1/4\pi^2,$ we note that $||u_s||_{\infty}$ is monotone increasing with $a$, whilst for $0<D<1/12\pi^2,$ we determine that now $||u_s||_\infty$ has a single maximum, at $a=\sqrt{3}\pi\sqrt{D}$, with value $||u_s||_\infty=1/(3^{\frac{3}{2}}\sqrt{D}).$ In Figure~\ref{fig3} we graph $||u_s||_\infty$ against $a$ at a number of values of $D.$ It is anticipated that, at fixed $0<D\le1/4\pi^2$, this branch of steady solutions will be smoothly continued into $a>1/2.$ This will be addressed at a later stage. For the present we now return to considering (DIVP) in detail.

 With equation (\ref{eqn1.1}) now having the form given above, we can, in fact, directly construct the solution to (DIVP) as follows. First we write the solution $u:\bar{D}_{\infty}\to \mathbb{R}$ to (DIVP), via Fourier's Theorem, as,
 \begin{equation}
    u(x,t) =  \sum_{n=1}^{\infty}{A_{n}(t)\sin\left(\frac{n\pi x}{a}\right)}~~\forall~~(x,t)\in \bar{D}_{\infty}, \label{eqn2.10}
 \end{equation}
 which satisfies the boundary conditions (\ref{eqn1.6}), and, for each $n\in\mathbb{N}$, $A_n:[0,\infty)\to \mathbb{R}$ is to be determined so that (\ref{eqn2.10}) satisfies equation (\ref{eqn2.2}) and initial condition (\ref{eqn1.5}). For each $n\in \mathbb{N}$, it is readily established that this is achieved if and only if $A_n(t)$ satisfies the initial value problem,
 \begin{equation}
     A_n'(t) = A_n(t) \left(1-D\frac{n^2\pi^2}{a^2}-a\bar{u}(t)\right),~~t>0, \label{eqn2.11}
 \end{equation}
 subject to,
 \begin{equation}
     A_n(0) = A_n^0 = \frac{2}{a}\int_0^a{u_0(y)\sin\left(\frac{n\pi y}{a}\right)}dy, \label{eqn2.12}
 \end{equation}
 and we note that, following the earlier conditions on $u_0,$
 \begin{equation}
     A_1^0>0,  \label{eqn2.13}
 \end{equation}
 whilst,
 \begin{equation}
     \bar{u}(0) = \frac{1}{a}\int_0^a{u_0(y)}dy = \frac{2}{\pi}\sum_{r=1}^{\infty}{\frac{A_{2r-1}^{0}}{(2r-1)}} > 0. \label{eqn2.14}
 \end{equation}
 The solution to (\ref{eqn2.11}) and (\ref{eqn2.12}) is readily obtained as,
 \begin{equation}
     A_n(t) = A_n^0J(t)e^{-\frac{n^2\pi^2D}{a^2}t},~~t\ge0,  \label{eqn2.15}
 \end{equation}
 with,
 \begin{equation}
     J(t) = e^{(t - aI(t))}>0,~~t\ge0,  \label{eqn2.16}
 \end{equation}
 and,
 \begin{equation}
     I(t) = \int_0^t{\bar{u}(s)}ds,~~t\ge0. \label{eqn2.17}
 \end{equation}
 The solution to (DIVP) can now be written as, via (\ref{eqn2.10}) and (\ref{eqn2.15}),
 \begin{equation}
     u(x,t) = J(t)\sum_{n=1}^{\infty}A_n^0e^{-\frac{n^2\pi^2D}{a^2}t}\sin\left(\frac{n\pi x}{a}\right),~~(x,t)\in \bar{D}_{\infty}. \label{eqn2.18}
 \end{equation}
 We note, via differentiating (\ref{eqn2.16}), that we may rearrange to obtain,
 \begin{equation}
     \bar{u}(t) = \frac{1}{a}\left(1 - \frac{J'(t)}{J(t)}\right),~~t\ge0,   \label{eqn2.19}  
 \end{equation}
 and, via (\ref{eqn2.16}), that,
 \begin{equation}
     J(0) = 1. \label{eqn2.20}
 \end{equation}
 Now, on integrating (\ref{eqn2.18}) across the spatial interval $[0,a]$, we obtain,
 \begin{equation}
     \bar{u}(t) = J(t)G(t),~~t>0, \label{eqn2.21}
 \end{equation}
 with $G\in C([0,\infty)) \cap C^{\infty}((0,\infty))$ being the known function,
 \begin{equation}
     G(t) = \frac{2}{\pi} \sum_{n=1}^{\infty}\frac{A_{2n-1}^0}{(2n-1)}e^{-\frac{(2n-1)^2\pi^2D}{a^2}t},~~t\ge0. \label{eqn2.22}
 \end{equation}
 It follows from (\ref{eqn2.16}) and (\ref{eqn2.21}) (recalling that it has been established earlier that the solution to (DIVP) is strictly positive on $D_{\infty}$) that,
 \begin{equation}
     G(t)>0~~\forall~~t\ge0,   \label{eqn2.23}
 \end{equation}
 and
 \begin{equation}
     G(0)=\bar{u}(0). \label{eqn2.24}
 \end{equation}
 whilst,
 \begin{equation}
     G(t) \sim \frac{2}{\pi}A_1^0e^{-\frac{\pi^2D}{a^2}t}~~\text{as}~~t\to \infty.  \label{eqn2.24'}
 \end{equation}
 Direct substitution from (\ref{eqn2.21}) into (\ref{eqn2.19}) then determines that $J(t)$ is the solution to the nonlinear Bernoulli equation,
 \begin{equation}
     J'(t) - J(t) = -aJ^2(t)G(t),~~t>0,  \label{eqn2.25}
 \end{equation}
 subject to the initial condition (\ref{eqn2.20}). This problem has solution $J:[0,\infty)\to \mathbb{R}$ given by,
 \begin{equation}
     J(t) = \frac{e^t}{(1 + aM(t))},~~t\ge0,  \label{eqn2.26}
 \end{equation}
 where $M:[0,\infty)\to \mathbb{R}$ is given by,
 \begin{equation}
    M(t) = \int_0^t{e^sG(s)}ds~~\forall~~t\ge0, \label{eqn2.27} 
 \end{equation}
 with $M(t)>0$ for $t>0.$ The exact solution to (DIVP) is now complete, and given by (\ref{eqn2.18}) with (\ref{eqn2.26}) and (\ref{eqn2.27}) On using (\ref{eqn2.24'}), we observe that, when $0<D<a^2/\pi^2$,
 \begin{equation}
     M(t)  \sim \frac{2A_1^0}{\pi}(1-\frac{\pi^2D}{a^2})^{-1}e^{(1-\frac{\pi^2D}{a^2})t}~~\text{as}~~t\to \infty,  \label{eqn2.28}
 \end{equation}
 whilst, when $D>a^2/\pi^2$,
 \begin{equation}
     M(t)\to M_{\infty}=\int_0^{\infty}{e^sG(s)}ds~~\text{as}~~t\to \infty. \label{eqn2.29}
 \end{equation}
 The marginal case when $D=a^2/\pi^2$ has,
 \begin{equation}
     M(t) \sim \frac{2A_1^0}{\pi}t~~\text{as}~~t\to \infty.  \label{eqn2.30}
 \end{equation}
 It now follows from (\ref{eqn2.27}) and (\ref{eqn2.28}) - (\ref{eqn2.30}) that,
\begin{equation}
    J(t) \sim
    \begin{cases}
        \frac{\pi}{2A_1^0a}\left(1 - \frac{\pi^2D}{a^2}\right)e^{\frac{\pi^2D}{a^2}t},~~0<D<\frac{a^2}{\pi^2}, \\ 
        \frac{e^t}{(1+aM_{\infty})},~~D>\frac{a^2}{\pi^2}, \\ 
        \frac{\pi e^t}{2A_1^0at},~~D=\frac{a^2}{\pi^2},
        \end{cases}
        \label{eqn2.31}
\end{equation}
as $t\to \infty.$ Thus we have, via (\ref{eqn2.18}) and (\ref{eqn2.31}), that the solution to (DIVP) has, when $0<D<a^2/\pi^2,$
\begin{equation}
    u(x,t) = u_s(x) + O(E(t))~~\text{as}~~t\to \infty, \label{eqn2.32}
\end{equation}
uniformly for $x\in [0,a]$, with $E(t)$ being exponentially small in $t$ as $t\to \infty.$ Conversely, when $D>a^2/\pi^2,$ we have,
\begin{equation}
    u(x,t)\sim\frac{A_1^0}{(1+aM_{\infty})}e^{-\left(\frac{\pi^2D}{a^2}-1\right)t}\sin\left(\frac{\pi x}{a}\right)~~\text{as}~~t\to\infty,  \label{eqn2.33}
\end{equation}
uniformly for $x\in[0,a]$. In the marginal case, $D=a^2/\pi^2,$ we have uniform decay to zero, as in the second case above, but now the decay is slower, being algebraic rather than exponential in $t$, as $t\to \infty$, with, specifically,
\begin{equation}
    u(x,t)\sim\frac{\pi}{2at}\sin\left(\frac{\pi x}{a}\right)~~\text{as}~~t\to \infty,  \label{eqn2.34}
\end{equation}
uniformly for $x\in[0,a].$

The cases with $0<a\le1/2$ are now complete, and we move on to consider the situation when $a>1/2.$

\subsection{$a>1/2$}
In this subsection we consider the remaining cases, which have $0<D<a^2/\pi^2$ when now $a>1/2.$ Although the scope for direct analytical progress is much more limited in this case, considerable progress can be made in a number of significant asymptotic limits, and this, complimented and supported by evidence gained from  detailed numerical solutions, results in a largely complete overall picture. Details of the simple numerical methods used to find steady solutions, compute bifurcation diagrams, solve the evolution problem and calculate eigenvalues for the linear stability problem, which are based on central finite differences and the trapezium rule, are given in Appendix~\ref{app_num}.

To begin with, we consider the non-negative and nontrivial steady states which can exist satisfying (\ref{eqn1.1}) with now (\ref{eqn1.3}), (\ref{eqn1.4}), and again the Dirichlet boundary conditions (\ref{eqn1.6}). This enables us subsequently to draw conclusions concerning the large-$t$ structure of the solution to (DIVP). It is again convenient to consider the steady state structure at fixed $D$ (as an unfolding parameter) whilst varying $a$ (as a bifurcation parameter). This leads us naturally to consider a number of complementary cases.

\subsubsection{$D\ge1/4\pi^2\approx0.02533$ with $0<(a-\pi\sqrt{D})\ll 1$}
In this case it is again straightforward to establish that a steady state transcritical bifurcation from the trivial equilibrium state, as $a$ increases through $\pi\sqrt{D}$, generates a nontrivial, non-negative steady state, and a weakly nonlinear theory is readily developed, which gives this bifurcated steady state the following structure close to the bifurcation point, namely,
\begin{equation}
    u_s(x;a,D) \sim A_b(D)(a-\pi\sqrt{D})\sin\left(\frac{\pi x}{a}\right)~~\text{as}~~a\to (\pi\sqrt{D})^+,  \label{eqn2.35}
\end{equation}
uniformly for $x\in[0,a]$. Here $A_b(D)$ is a strictly positive $O(1)$ constant, independent of $a,$ and, for brevity, we omit the standard details.

\subsubsection{$0<D<1/4\pi^2$ with $0<(a-1/2)\ll1$}
In this case we consider whether the unique steady state we have identified in subsection 2.1, which exists for $\pi\sqrt{D}<a\le1/2,$ can be continued into $a>1/2$. It is readily established that this continuation is certainly possible, and uniquely so, at least locally, when $a$ is sufficiently close to 1/2. The continued steady state has the local asymptotic structure,
\begin{equation}
    u_s(x;a,D) \sim \pi(1-4\pi^2D)\sin\left(\frac{\pi x}{a}\right)~~\text{as}~~ a\to \frac{1}{2}^+,  \label{eqn2.36}
\end{equation}
uniformly for $x\in [0,a]$.
In the next case we consider steady states when $D$ is large, with $a=O(\sqrt{D})$ as $D\to \infty$.

\subsubsection{$D\gg1$ with $a=O(\sqrt{D})$}
We write,
\begin{equation}
    a = \sqrt{D}\bar{a},   \label{eqn2.37}
\end{equation}
with $\bar{a}=O(1)^+$ as $D\to \infty.$ We anticipate that steady states have $u_s=O(1)$ as $D\to \infty.$ Thus a nontrivial balance in equation (\ref{eqn1.1}), as $D\to \infty$, requires us to introduce the scaled coordinate,
\begin{equation}
    \bar{x} = x/\sqrt{D},  \label{eqn2.38}
\end{equation}
with now $\bar{x}\in [0,\bar{a}]$. For a steady state, the nonlocal term becomes (except in passive boundary layers, where $\bar{x}=O(1/\sqrt{D})^+ \text{or}~ \bar{x}= \bar{a}-O(1/\sqrt{D})^+$),
\begin{equation}
 \sqrt{D}\int_{\bar{x}-\frac{1}{2\sqrt{D}}}^{\bar{x}+\frac{1}{2\sqrt{D}}}{u_s(\bar{y})}d\bar{y} = u_s(\bar{x}) + O(D^{-1})~~\text{as}~~D\to \infty,
 \label{eqn2.39}
\end{equation}
with $\bar{x}\in (O(1/\sqrt{D})^+, \bar{a}-O(1/\sqrt{D})^+).$ A steady state $u_s:[0,\bar{a}]\to \mathbb{R}$ is now expanded in the form,
\begin{equation}
    u_s(\bar{x},\bar{a};D) = U_s(\bar{x},\bar{a}) + O(D^{-1})~~as~~D\to \infty,   \label{eqn2.40}
\end{equation}
with $\bar{x}\in [0,\bar{a}]$. On substitution from expansion (\ref{eqn2.40}) into equation (\ref{eqn1.1}) (when written in terms of the coordinate $\bar{x}$), we obtain the following local, nonlinear ODE for $U_s$, namely,
\begin{equation}
    U_s'' + U_s(1-U_s) = 0,~~\bar{x}\in (0,\bar{a}),   \label{eqn2.41}
\end{equation}
where $'=d/d\bar{x}$. We note that this is precisely the Dirichlet problem for steady states on the interval $[0,\bar{a}]$ for the \emph{local} Fisher-KPP equation (see, for example, \cite{PCF}). This is to be solved subject to the boundary and non-negativity conditions,
\begin{equation}
    U_s(0,\bar{a})=U_s(\bar{a},\bar{a})=0~~,~~U_s(\bar{x},\bar{a})\ge0~~\forall~~x\in[0,a].   \label{eqn2.42}
\end{equation}
It is straightforward to consider solutions to this nonlinear boundary value problem, which we label as (NBVP), by recasting it in the $(U_s,U_s')$ phase plane. We omit details, and simply record the conclusions. Firstly, (NBVP) has no nontrivial solutions when $0<\bar{a}\le \pi$. However, consistent with the earlier theory, a transcritical bifurcation from the trivial solution occurs in (NBVP) as $\bar{a}$ increases through the bifurcation value $\bar{a}=\pi$ (this being, in term of $a$, the value $a=\pi\sqrt{D},$ consistent with the earlier theory). This bifurcation produces a unique nontrivial solution to (NBVP) in $\bar{a}>\pi.$ We denote this solution by $U_s=\nu(\bar{x},\bar{a})$ with $\bar{x}\in [0,\bar{a}]$, and it has the following key properties,
    \begin{enumerate}
    \item $\nu(\bar{x},\bar{a})>0~~\forall~~\bar{x}\in (0,\bar{a}),$
    \item $\nu(\bar{x},\bar{a})$ is symmetric about $\bar{x}=\frac{1}{2}\bar{a},$
    \item $\nu(\bar{x},\bar{a})$ has a single turning point, which is a maximum, at $\bar{x}=\frac{1}{2}\bar{a}$,
    \item with $A(\bar{a})=\text{sup}_{\bar{x}\in [0,a]}(\nu(\bar{x},\bar{a}))$, then $A(\bar{a})$ is monotone increasing with $\bar{a}>\pi$, and $A(\bar{a})\to 0~ \text{as}~ \bar{a}\to \pi$, whilst $A(\bar{a})\to 1$ as $\bar{a}\to \infty.$
\end{enumerate}
It is also straightforward to establish that,
\begin{equation}
    \nu(\bar{x},\bar{a}) = \frac{3}{4}(\bar{a}-\pi)\sin\left(\frac{\pi\bar{x}}{\bar{a}}\right)+O((\bar{a}-\pi)^2)~~\text{as}~~\bar{a}\to \pi,
\end{equation}
uniformly for $\bar{x}\in [0,\bar{a}],$ whilst,
\begin{equation}
    \nu(\bar{x},\bar{a}) = 
    \begin{cases}
     1 - \frac{3}{2}\text{sech}^2 \left\{\frac{1}{2}\left(\bar{x}+\ln{(2+\sqrt{3})}\right)\right\}+o(1),~~\bar{x}=O(1)^+,\\
    1 - \frac{6}{2+\sqrt{3}}\left( e^{-\bar{x}} + e^{-(\bar{a}-\bar{x})}\right)~+~o(e^{-\bar{x}},~e^{-(\bar{a}-\bar{x})}),~~\bar{x}=O(\bar{a})^+,\\
    1 - \frac{3}{2}\text{sech}^2\left\{\frac{1}{2}\left((\bar{a}-\bar{x})+\ln{(2+\sqrt{3})}\right)\right\}+o(1),~~\bar{x}=\bar{a}-O(1)^+,
    \end{cases}
    \label{eqn2.44}
\end{equation}
as $\bar{a}\to \infty.$ This concludes the analysis for $D$ large with $a=O(\sqrt{D}),$ which has established that there are no non-negative nontrivial steady states for $0<a\le \pi\sqrt{D}$, but a unique such steady state at each $a>\pi\sqrt{D}.$\\\\
The next case has $D$ of $O(1)$ with $a$ large.

\subsubsection{$0<D\le O(1)$ with $a\gg1$}
\begin{figure}
\begin{center}
\includegraphics[width= 0.8\textwidth]{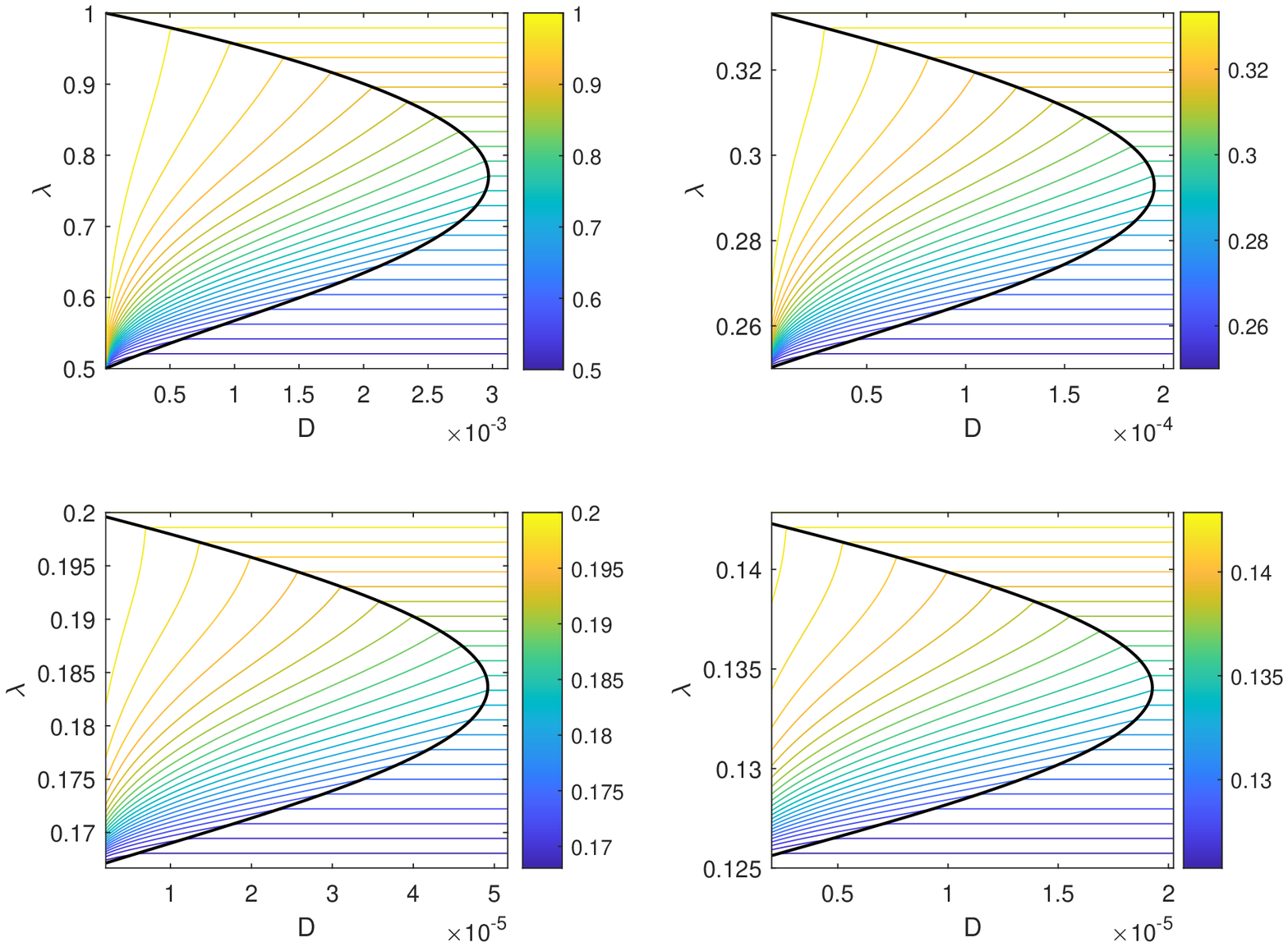}
\caption{The $L^1$ norm, $A(\lambda,D)$, of the nontrivial, non-negative, periodic steady states identified in (NB1). The first four tongue regions are shown in the $(D,\lambda)$ plane.}\label{heatmap}
\end{center}
\end{figure}

\begin{figure}
\begin{center}
\includegraphics[width=0.8\textwidth]{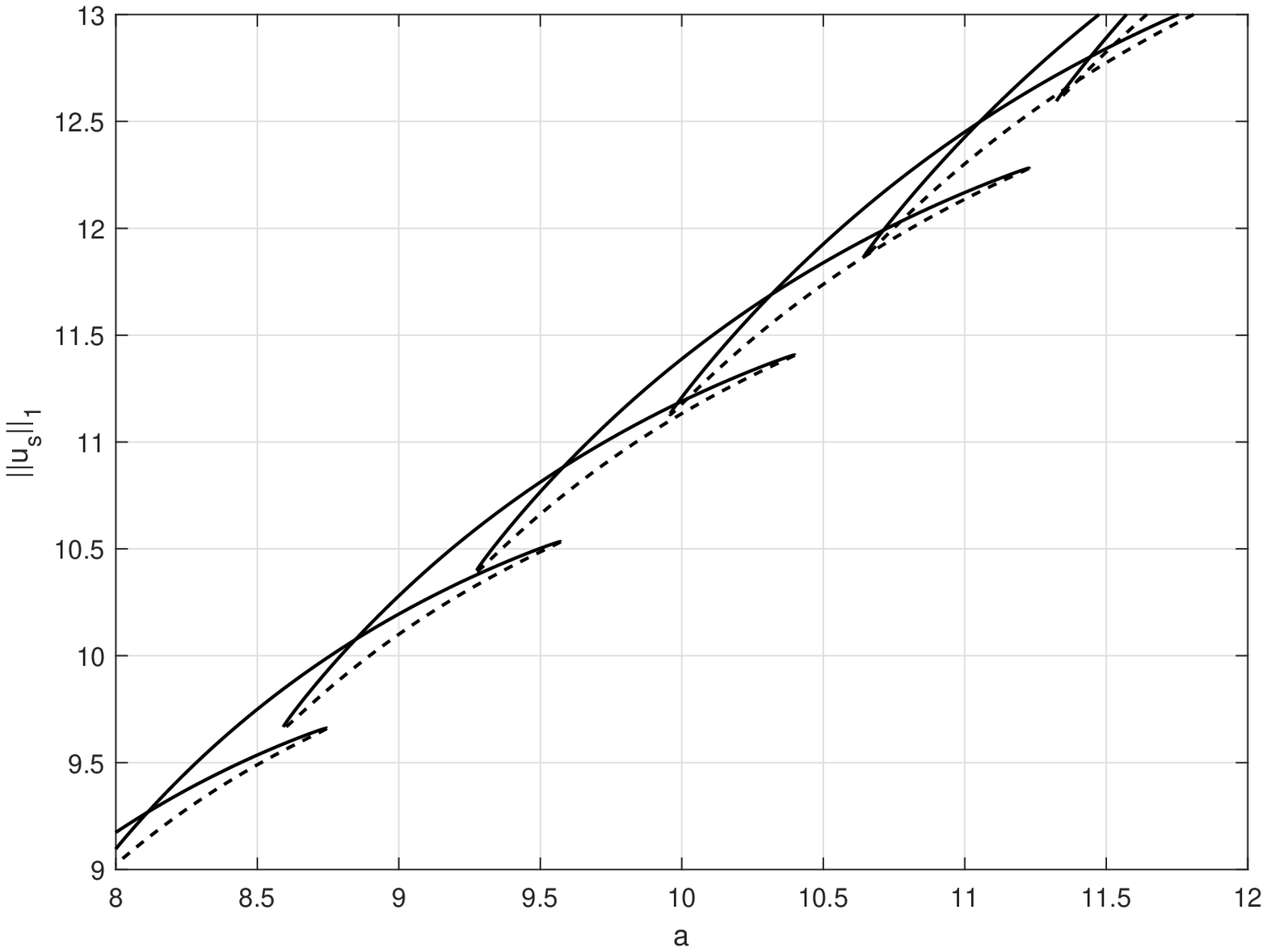}
\caption{Part of the $(a,||u_s||_1)$ bifurcation diagram, when $D = 2 \times 10^{-3}$. Unstable branches of steady states are indicated by broken lines.}\label{D2em3_bif}
\end{center}
\end{figure}

In this case, since the domain length $a$ is large, whilst $D$ remains of O(1), we anticipate that there will be symmetric boundary layers when $x=O(1)^+$ and $x=a-O(1)^+$ as $a\to \infty$, at either end of a bulk region separating these two boundary layers. We expect that any nontrivial steady state will have $u_s=O(1)$ in the bulk region and in the boundary layers, as $a\to \infty$. It is instructive to begin in the boundary layers, and as these are symmetric we need only consider the boundary layer at the left end of the domain. In this boundary layer we write,
\begin{equation}
    u_s(x,D,a) \sim U_b(x,D)~~\text{as}~~a\to \infty,  \label{eqn2.45}
\end{equation}
with $x=O(1)^+.$ The leading order problem is then,
\begin{equation}
    DU_b'' + U_b\left(1 - \int_{\alpha(x)}^{x+\frac{1}{2}}{U_b(y,D)}dy\right) = 0,~~x=O(1)^+,   \label{eqn2.46}
\end{equation}
 with $U_b$ being non-negative, and subject to the boundary conditions,
\begin{equation}
    U_b(0,D)=0~~\text{and}~~U_b(x,D) ~\text{bounded as}~~x\to \infty.   \label{eqn2.48}
\end{equation}
A numerical and analytical investigation of this nonlocal boundary value problem, which we label as (SIBP), establishes the following key properties :

\begin{enumerate}
\item For each $D\ge \Delta_1\approx 0.00297$ (see (NB1), section 4) (SIBP) has a unique solution. The solution has $U_b(x,D)\to 1$, exponentially in $x$,  as $x\to \infty$. In addition, $U_b(x,D)$ is monotone when $D\ge \Delta_m$, but has harmonic exponential decay when $\Delta_1<D<\Delta_m$, with the exponential decay rate reducing to zero as $D\to \Delta_1$ . Here $\Delta_m=2 l^{-3}\sinh\left(\frac{1}{2}l\right)$, with $l \approx 5.9694$ being the unique positive solution to the algebraic equation $6\tanh\left(\frac{1}{2}l\right) = l$. This determines $\Delta_m \approx 0.0928$.

\item For each $0<D<\Delta_1$, (SIBP) now has a one-parameter family of solutions, parameterised by a characteristic wavelength $\lambda$, with $\lambda\in (\lambda_1^-(D),\lambda_1^+(D))$, where the notation here is as introduced in section 4 and section 5 of (NB1). Specifically, each of these solutions is now asymptotic to a steady non-negative, periodic solution of (\ref{eqn1.1}), oscillating  about unity, and with wavelength $\lambda$, as $x\to \infty$. In particular
\begin{equation}
  U_b(x,D)\sim F_p(x - x_b(\lambda,D),\lambda,D)~~\text{as}~~x\to \infty.   \label{eqn2.49}  
\end{equation}
Here $F_p(y,\lambda,D)$, for $y\in [0,\lambda)$ and $(\lambda,D)\in \Omega_1$, is that nontrivial, non-negative periodic steady state, with wavelength $\lambda$, and oscillating about unity, as introduced in section 5 of (NB1), with the translation fixed so that $F_p(0,\lambda,D)=1~~\text{and}~~ F_p'(0,\lambda,D)>0$. In addition, $x_b(\lambda,D)$ is a shift of origin which is determined in the solution of (SIBP).
\end{enumerate}
This completes the structure in the left boundary layer. In the right boundary layer, we have, symmetrically,
\begin{equation}
u_s(x,D,a)\sim U_b(a-x,D)~~\text{as}~~a\to \infty,   \label{eqn2.50}
\end{equation}
with $x=a-O(1)^+$. We must now connecct the two boundary layers through the bulk region, when $x\in (O(1)^+, a-O(1)^+)$. The situation for $D\ge \Delta_1$ is straightforward. In the bulk region we must simply take,
\begin{equation}
  u_s(x,D,a) \sim 1 ~~\text{as}~~a\to \infty, \label{eqnm}  
\end{equation}
for $x\in (O(1)^+,a-O(1)^+)$, and the structure is complete. For $0<D<\Delta_1$ we begin by recalling from (NB1) the following  general properties of $F_p(x,\lambda,D)$ with $x\in \mathbb{R}$:
\begin{itemize}
    \item $F_p(x,\lambda,D)$ has exactly one maximum point and one minimum point, and no other stationary points, over one wavelength $\lambda$
    
    \item Consecutive maximum and minimum points are separated by a distance $\frac{1}{2}\lambda$.
    
    \item $F_p(x,\lambda,D)$ is an even function \emph{about} each of its stationary points.
\end{itemize}
Now let $x=x_M(\lambda,D)$ be the location of the first maximum point of $F_p(x,\lambda,D)$ on the positive $x$-axis. Using the above properties, it is readily established, for each $r\in \mathbb{N}$, $D\in (0,\Delta_1)$ and $\lambda\in (\lambda_1^+(D),\lambda_1^-(D))$, that the function $F_p(x-x_b(\lambda,D),\lambda,D)$ has a stationary point at the positive location
\begin{equation}
x = x_r(\lambda,D) = x_b(\lambda,D) + x_M(\lambda,D) +\frac{1}{2}(r-1)\lambda. \label{eqn2.51}
\end{equation}
Therefore, following the final property above, we may conclude that $F_p(x-x_b(\lambda,D),\lambda,D)$:
\begin{itemize}
\item is an even function of $x$ \emph{about} $x=x_r(\lambda,D)$,
\item has a maximum or minimum point at this location when $r$ is odd or even respectively,
\item has exactly $r$ maximum points for $x\in [0,2x_r(\lambda,D)]$.
\end{itemize}
We can now return to the bulk region. First fix $D\in (0,\Delta_1)$, and choose $r\in \mathbb{N}$ with $r$ large. Then for each $\lambda\in (\lambda_1^-(D),\lambda_1^+(D))$, take,
\begin{equation}
a = 2x_r(\lambda,D) = \lambda(r-1) + O(1)~~\text{as}~~r\to \infty.\label{eqn2.52}
\end{equation}
For this large value of $a$, as $r\to \infty$, we then take,
\begin{equation}
u_s(x,D,a) \sim F_p(x-x_b(\lambda,D),\lambda,D) \label{eqn2.53}
\end{equation}
for $x\in (O(1)^+,a-O(1)^+)$. We observe, via (\ref{eqn2.52}) and the first comment above, that (\ref{eqn2.53}) matches with (\ref{eqn2.45}) and (\ref{eqn2.50}) when $x=O(1)^+$ and $x=a-O(1)^+$ respectively. Thus, to summarise, for each fixed $0<D<\Delta_1$ we conclude, for each large natural number $r$, and then each
\begin{equation}
a \in (\lambda_1^-(D)(r-1)+O(1)^+, \lambda_1^+(D)(r-1)+O(1)^+) \equiv I_r(D),  \label{eqn2.54}
\end{equation}
that a nontrivial, non-negative steady state exists. Away from boundary layers at the ends of the long domain, this solution is periodic about unity, with wavelength (for $a=O(r)$ as given in (\ref{eqn2.54}))
\begin{equation}
\lambda \sim a(r-1)^{-1},   \label{eqn2.55}
\end{equation}
and having $r$ peak points across the domain. To characterise each such steady state it is instructive to measure its norm in $L^1$. Firstly, for $0<D<\Delta_1$ and $\lambda\in (\lambda_1^-(D),\lambda_1^+(D))$, we write, over an interval of one wavelength $\lambda$,
\begin{equation}
A(\lambda,D) = || F_p(\cdot,\lambda,D)||_1 \equiv \int_0^\lambda{F_p(y,\lambda,D)}dy,  \label{eqn2.56} 
\end{equation}
and numerically determined contour plots of $A(\lambda,D)$ in the $(D,\lambda)$ plane, for each of the first four tongues where the periodic steady states exists (see (NM1), section 5), are shown in Figure~\ref{heatmap}. Outside the tongues, $A$ varies linearly with $\lambda$ as we have plotted the result for the equilibrium solution, $u_e=1$. Within the tongues, the deviation of $A$ from this linear behaviour becomes stronger as $D$ decreases, becoming independent of $\lambda$ as $D \to 0$. The qualitative behaviour of $A(\lambda,D)$ is very similar in each tongue.
It is now straightforward to determine, via construction and using $A(\lambda,D)$ in the first tongue, that for fixed $D$, and given large $r$, then,
\begin{equation}
||u_s(\cdot,D,a)||_1 \sim rA(a(r-1)^{-1},D)~~\forall~~a\in ((r-1)\lambda_1^-(D),(r-1)\lambda_1^+(D)).  \label{eqn2.57}
\end{equation}
With $D$ fixed and $r$ large, each interval $I_r(D)$, containing the branch of $r$-peak steady states, has an echelon overlap with its predecessor $I_{r-1}(D)$, containing the branch of $(r-1)$-peak steady states, and this structure continues as $r$ increases, with the number of overlaps, with previous echelons, increasing linearly with increasing $r$. This can be represented on the $(a,||u_s(\cdot,D,a)||_1)$ bifurcation plane. For each large $r$, the associated $r$-peak branch of steady states has a graph given by (\ref{eqn2.57}) on this plane. The graphs of these branches form a lifting, and multiply  overlapping, echelon of separate segments. At a given large $a$, the ovelapping echelons are associated with those natural numbers $r_m(a,D)\le r\le r_M(a,D)$, with
\begin{equation}
 r_m(a,D) \sim \frac{a}{\lambda_1^+(D)}  + 1,
\end{equation}
\begin{equation}
    r_M(a,D) \sim \frac{a}{\lambda_1^-(D)} + 1,
\end{equation}
as $a\to \infty$. Each segment has a base length increasing linearly with $r$, a slope of O(1) and a consequent $O(r)$ increase in norm on traversing the segment. In addition each segment is lifted by an $O(1)$ value each time $r$ is increased by unity. Careful numerical solution of the full boundary value problem, using pseudo-arclength continuation in the $(a, ||u_s||_1)$-plane, reveals that each two consecutive, overlapping, branches in the echelon are continuously connected by a hanging, elongated loop with two facing saddle-node bifurcations, as illustrated in Figure~\ref{D2em3_bif}. On the lower part of this loop 
\begin{sidewaysfigure}
\begin{center}
\vspace{12.5cm}
\includegraphics[width=\columnwidth]{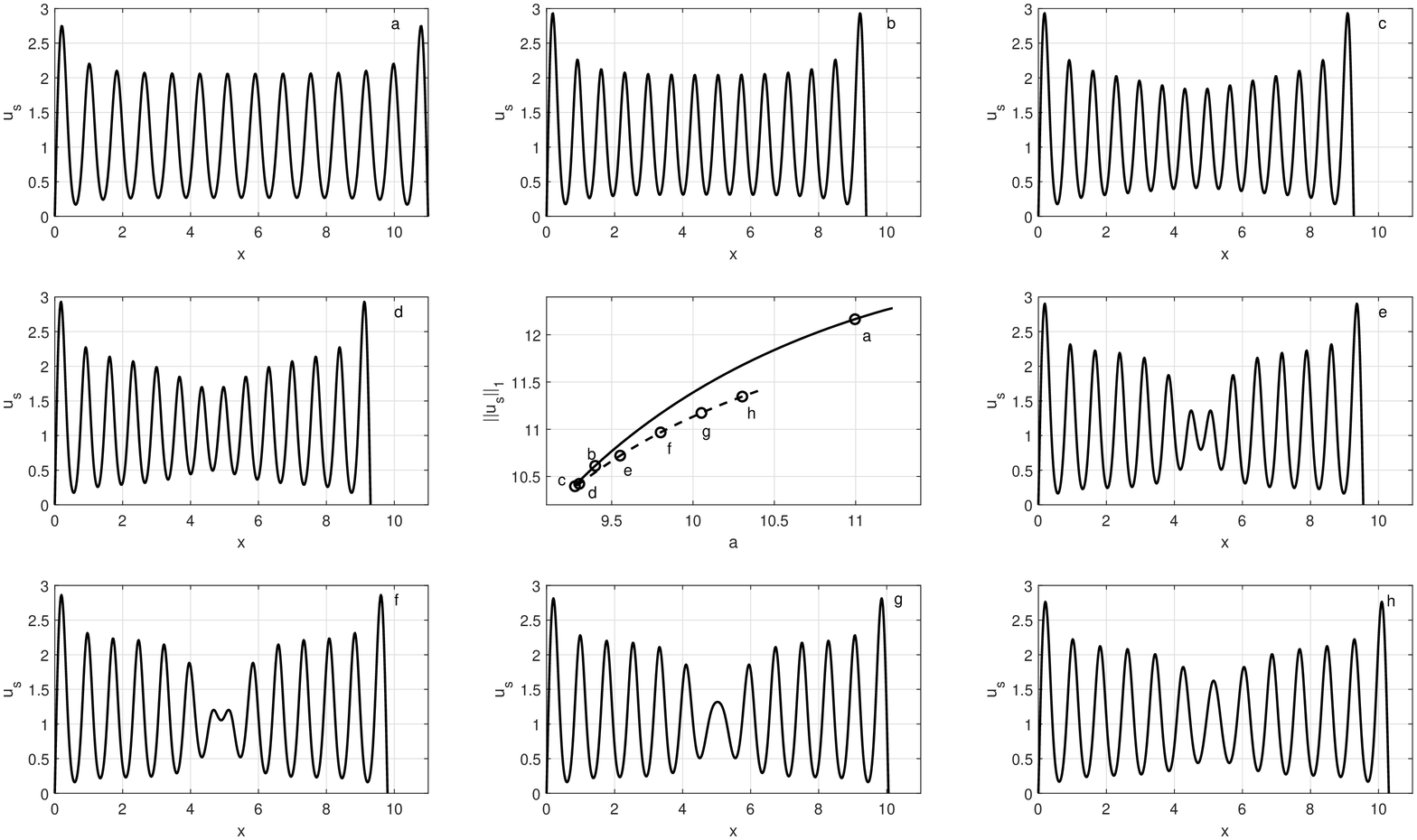}
\caption{The profile of the steady state $u_s$, with $D=2\times 10^{-3}$, as the $r=14$ to $r=13$ echelon of the bifurcation diagram is traversed. The central panel shows the relevant portion of the bifurcation diagram. Note the different $u_s$-axis range in the other panels, individually labelled (a) to (h), and identified on the central panel. Also only the  steady states (a) to (c), on the upper echelon, are stable, with the steady states (d) to (h), on the connecting loop, each being unstable. }\label{D2em3_sols}
\end{center}
\end{sidewaysfigure}
the associated steady state develops an additional peak through a rapid process which occurs, for $r$ large,  in the neighbourhood of the stationary point at  $x=\frac{1}{2}a$, in the bulk region, and allows this central stationary point to flip through itself. There is then a rapid transition to the associated multi-peak steady state as each respective saddle-node bifurcation is traversed. A typical case is shown in Figure~\ref{D2em3_sols}.

The above structure continues to hold as $D$ becomes small, with $a$ large. However, when $D$ is small we can take advantage on the approximations to $F_p(x,\lambda,D)$ as detailed in section 5 of (NB1). In particular we have, as $D\to 0$,
\begin{equation}
(\lambda_1^-(D),\lambda_1^+(D))\to \left(\frac{1}{2},1\right),  \label{eqn2.58}
\end{equation}
\begin{equation}
 x_b(\lambda,D) + x_M(\lambda,D) \to \frac{1}{2}\left(\lambda-\frac{1}{2}\right)~~\forall~~\lambda\in \left(\frac{1}{2}+O(\sqrt{D}),1-O(D)\right),     \label{eqn2.59}
\end{equation}
and
\begin{equation}
A(\lambda,D)\to 1~~\forall~~\lambda\in \left(\frac{1}{2}+O(\sqrt{D}),1-O(D)\right).   \label{eqn2.60}
\end{equation}
This completes the examination of the steady state structure to (DIVP) when $a$ is large. \\\\We now continue by considering the steady state structure when $D$ is small and $a>1/2 + o(1)^+$ as $D\to 0$.

\subsubsection{$0<D\ll1$ with $a>1/2 + o(1)$ as $D\to 0$}
\begin{figure}
\begin{center}
\includegraphics[width=\textwidth]{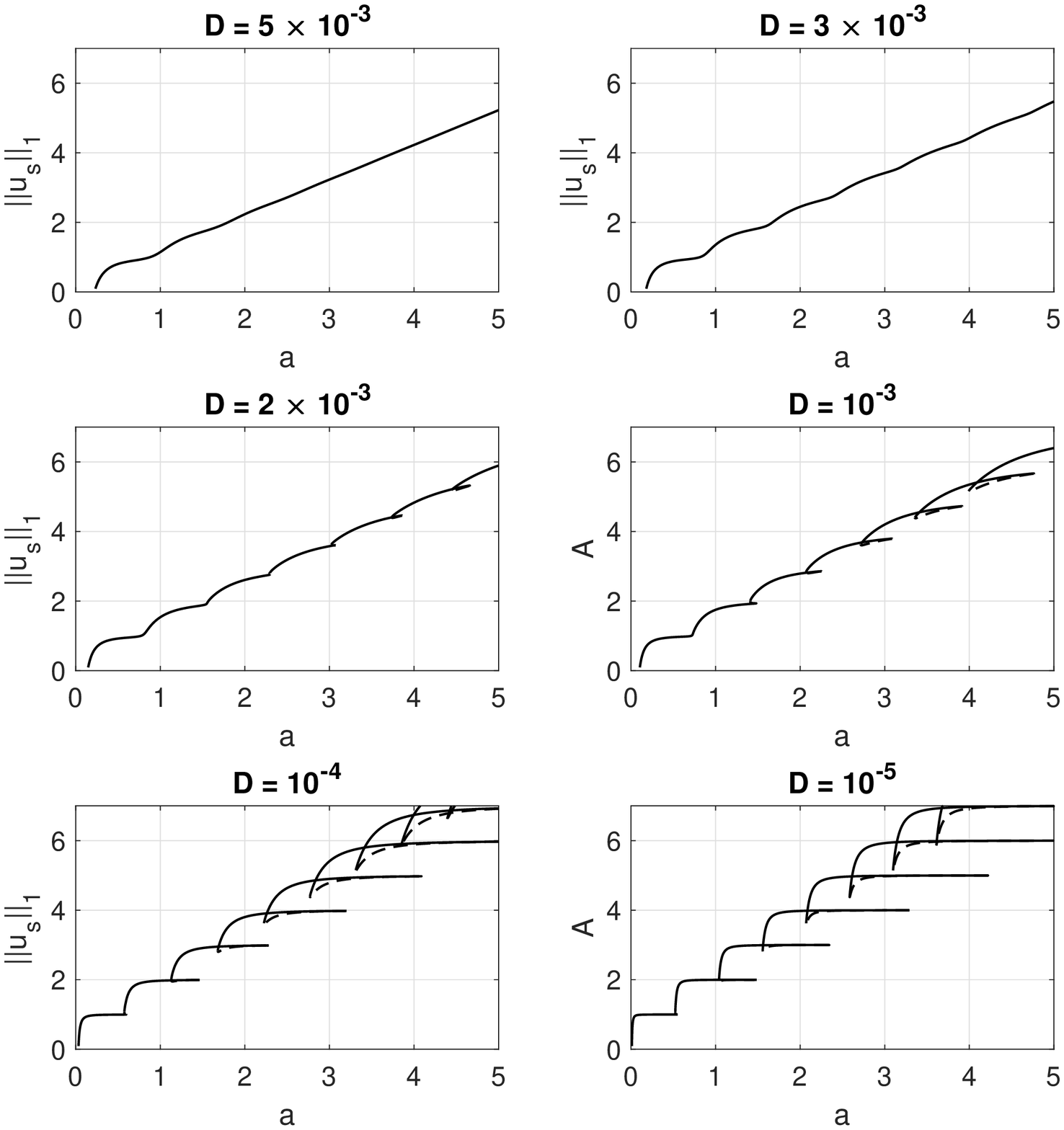}
\caption{The $(a,||u_s||_1)$ bifurcation diagram for various values of $D$. Unstable branches are shown as broken lines.}\label{fig2}
\end{center}
\end{figure}

\begin{sidewaysfigure}
\begin{center}
\vspace{12.5cm}
\includegraphics[width=\columnwidth]{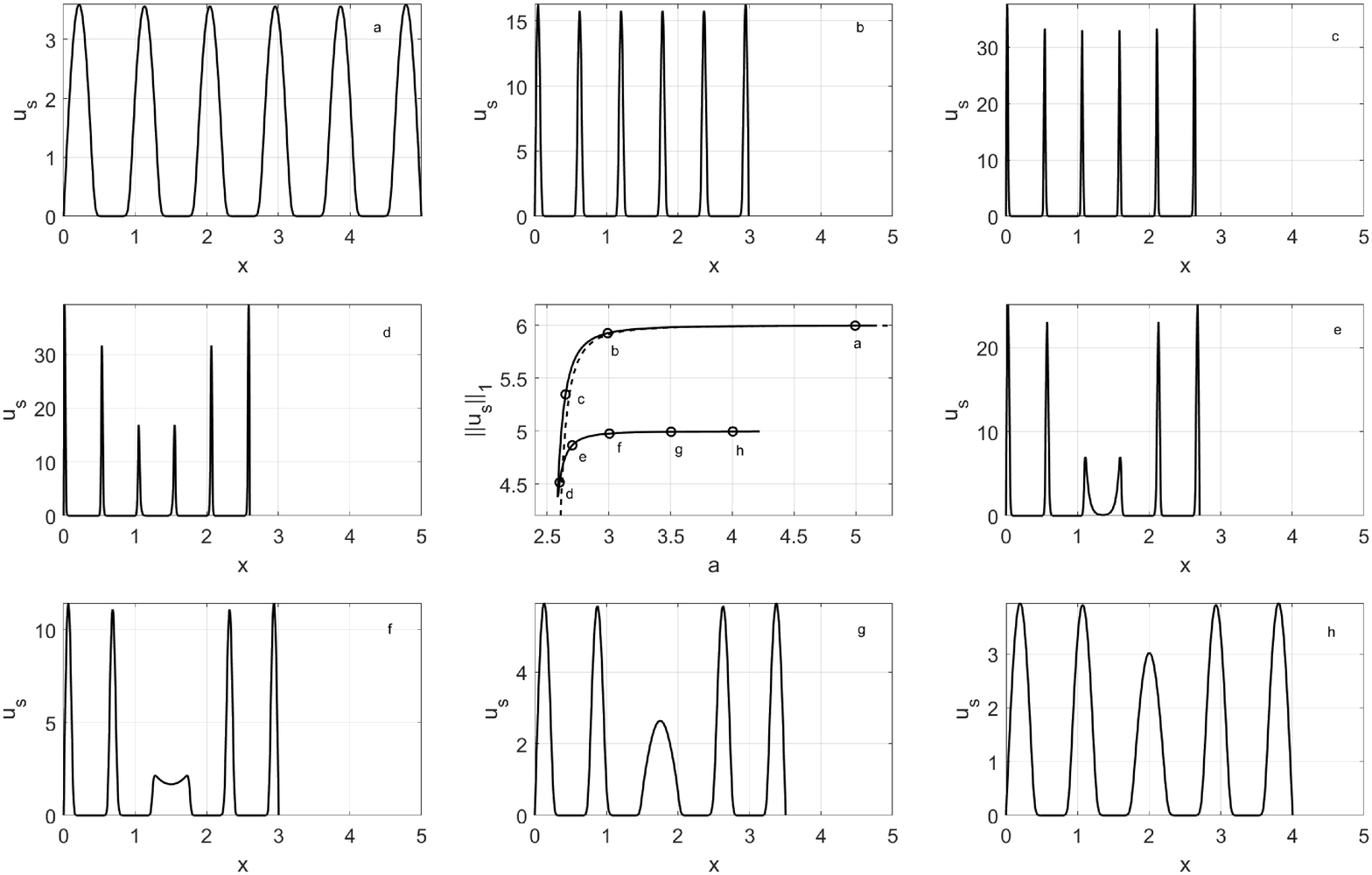}
\caption{The profile of the steady state $u_s$ as the $r=6$ to $r=5$ echelon of the bifurcation diagram is traversed. The central panel shows the relevant portion of the bifurcation diagram, with the broken line given by (\ref{eqn2.66}). Note the different $u_s$-axis limits in the other panels, individually labelled (a) to (h), and identified on the central panel. The most dramatic changes in functional form occur along the upper branch of the bifurcation curve, where the peaks become taller and thinner, and then at the left hand side of the bifurcation curve, where the number of peaks changes by one. Note that only solutions (a) to (c) are stable. }\label{fig_U}
\end{center}
\end{sidewaysfigure}

\begin{sidewaysfigure}
\begin{center}
\vspace{12.5cm}
\includegraphics[width=\columnwidth]{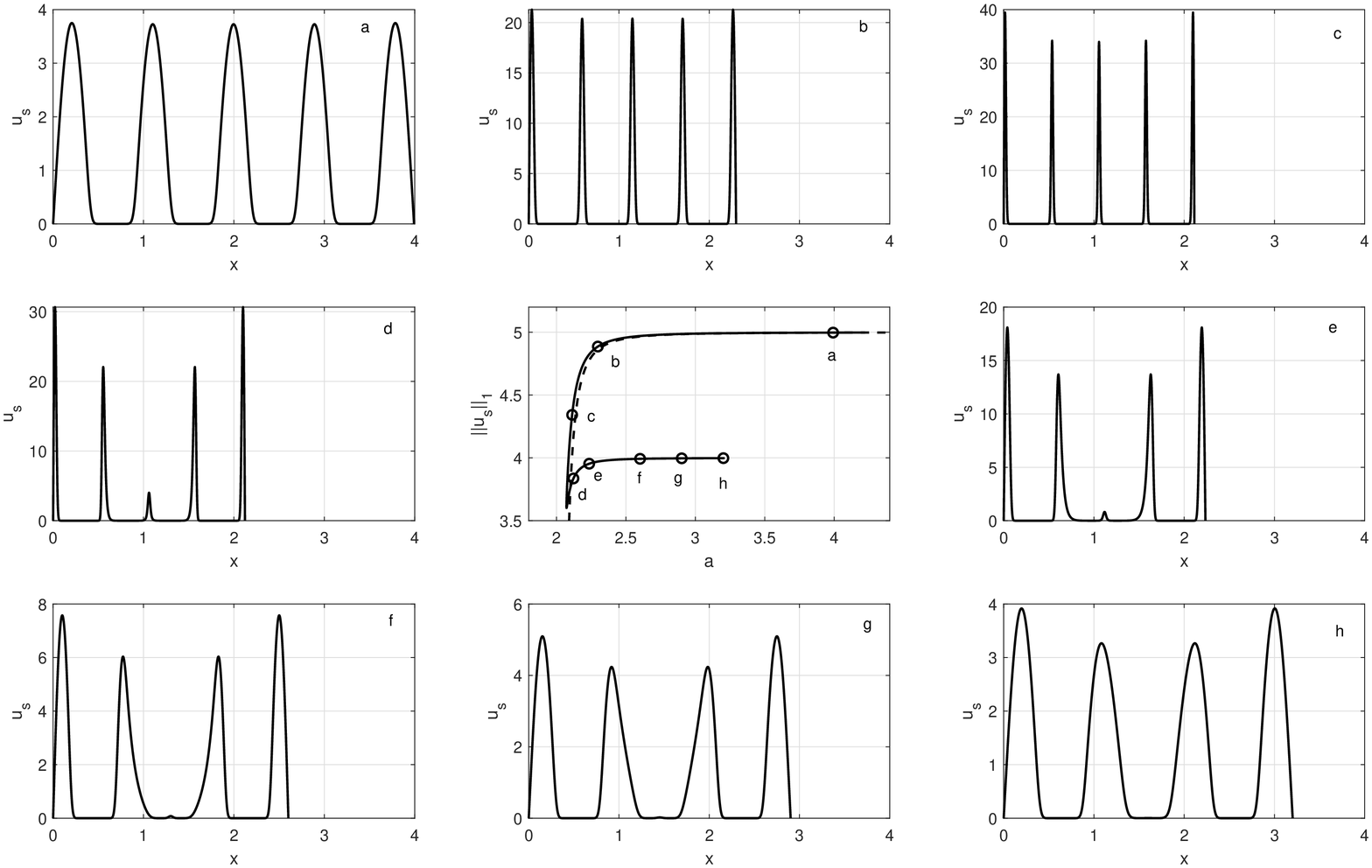}
\caption{The profile of the steady state $u_s$ as the $r=5$ to $r=4$ echelon of the bifurcation diagram is traversed. The central panel shows the relevant portion of the bifurcation diagram, with the broken line given by (\ref{eqn2.66}). The steady state in (a) of  this Figure is close to the steady state in (h) of the previous Figure~\ref{fig_U} (see Figure~\ref{fig2}, bottom right hand panel) and is qualitatively very similar. Note the different $u_s$-axis limits in the other panels, individually labelled (a) to (h), and identified on the central panel. Note that only solutions (a) to (c) are stable.}\label{fig_U2}
\end{center}
\end{sidewaysfigure}

In this case, we take advantage of the detailed theory developed in subsection $5.1$ of (NB1) for the periodic states $F_p(x,\lambda,D)$ with $\lambda\in (\lambda_1^-(D),\lambda_1^+(D))$, using the terminology laid out in (NB1) in subsection $5.1$. Summarising from (NB1), the family of nontrivial and non-negative periodic states are given by (at leading order as $D\to 0$, and noting that a passive edge region, of thickness $O(D^{\frac{1}{4}})$, is present at the end of the supported leading order form for $F_p$, which provides local smoothing to the gradient discontinuity in the leading order form), for each wavelength $\lambda\in (1/2 + O(D^{\frac{1}{2}})^+, 1-O(D)^+)$,
\begin{equation}
F_p\left(x - \frac{1}{2}\left(\lambda-\frac{1}{2}\right),\lambda,D\right) \sim
\begin{cases}
\frac{\pi}{(2\lambda-1)}\sin\left(\frac{\pi x}{\lambda-\frac{1}{2}}\right),~~~ 0\le x \le \lambda-\frac{1}{2}, \\
0,~~~\lambda-\frac{1}{2} < x \le \lambda,
\end{cases}
\label{eqn2.61}
\end{equation}
as $D\to 0$ with $x\in [0,\lambda]$. We also recall that,
\begin{equation}
A(\lambda,D) = 1 - \frac{\pi^2}{(\lambda-\frac{1}{2})^2}D + o(D^{\frac{5}{4}})~~\text{as}~~D\to 0.  \label{eqn2.62}
\end{equation}
with $\lambda\in (1/2+O(D^{\frac{1}{2}})^+,1-O(D)^+)$. Now, for each $r=2,3,4...$, we can construct from (\ref{eqn2.61}) an $r$-peak steady state at each domain size $a$ with
\begin{equation}
a \in \left(\frac{1}{2}(r-1) + O(D^{\frac{1}{2}})^+,r-\frac{1}{2}-O(D)^+\right),   \label{eqn2.63}
\end{equation}
which is given by
\begin{equation}
u_s(x,D,a) \sim F_p\left(x - \frac{1}{2}\left(\lambda(a,r)-\frac{1}{2}\right),\left(a+\frac{1}{2}\right)r^{-1},D\right)~~\text{as}~~D\to 0,   \label{eqn2.64}
\end{equation}
for $x\in [0,a]$, and $F_p$ approximated as in (\ref{eqn2.61}), with
\begin{equation}
\lambda(a,r)=\left(a+\frac{1}{2}\right)r^{-1}.  \label{eqn2.65}
\end{equation}
The steady state is uniformly approximated by  (\ref{eqn2.64}), except in passive, thin,  boundary layers at each end of the domain, when $x=O(D^{\frac{1}{2}})^+$ or $x=a-O(D^{\frac{1}{2}})^+$, respectively, and in which $u_s=O(D^{\frac{1}{2}})$. In the bifurcation plane each branch of $r$-peak steady states lies on the curve,
\begin{equation}
||u_s(\cdot,D,a)||_1 = r\left(1 - \frac{\pi^2}{\left(\lambda(a,r)-\frac{1}{2}\right)^2}D + o(D^{\frac{5}{4}})\right)~~\text{as}~~D\to 0   \label{eqn2.66}
\end{equation}
for $a \in \left(\frac{1}{2}(r-1) + O(D^{\frac{1}{2}})^+,r-\frac{1}{2}-O(D)^+\right)$, as illustrated in the central panels of Figures~\ref{fig_U} and~\ref{fig_U2}, where this asymptotic expression can be seen to be in excellent agreement with the numerically-calculated bifurcation curve. These overlapping branches form an increasing echelon as $r$ is successively increased, and each successive echelon is raised by unity, at leading order as $D\to 0$. Each pair of overlapping echelon branches is connected by a \emph{loop}, which consists of a very stiff subcritical saddle-node bifurcation at the head of each echelon, after which it then traces the echelon  backwards, finally increasing towards the tail of the echelon, passing through a less stiff supercritical saddle-node bifurcation and then increasing by approximately unity to connect with the next echelon. Whilst traversing the loop, from head to tail, the additional peak is generated by either the minimum or maximum point at $x=\frac{1}{2}a$ \emph{pulling through itself}, when $r$ is even or odd respectively. This structure has been verified using pseudo-arclength continuation, and a typical numerically determined bifurcation diagram, for $D=10^{-5}$, is shown in the final window of Figure~\ref{fig2}. 
In addition, for the echelon - loop structure consisting of the $r=5$ and $r=6$ echelons, and then the $r=4$ and $r=5$ echelons, together with their corresponding loop connections, we graph, in Figure~\ref{fig_U} and in Figure~\ref{fig_U2}, a sequence of the associated steady states as the echelon - loop - echelon structure is traversed. In both Figures we give profiles from the upper echelon, moving backwards along this echelon onto the loop and through the lower supercritical saddle-node bifurcation, when the steady state begins its transition from a $(r+1)$-peak to an $r$-peak form and then continue forwards along the loop towards the upper subcritical saddle-node bifurcation. We observe that the steady state profile is close to achieving its final $r$-peak form prior to going through this upper saddle-node point, and so, for brevity, we stop giving profiles in both Figures towards the upper end of the lower loop.
We note the excellent agreement with the above theory for $D$ small. We also note, from (\ref{eqn2.66}) above, that we should anticipate that the supercritical saddle-node bifurcation should take place in a region with $(a,||u_s||_1) = \left(\frac{1}{2}(r-1) -O(\sqrt{D}), r - O(1)\right)$, whilst the subcritical saddle-node bifurcation should take place in a tighter region with $(a,||u_s||_1) = \left(r-\frac{1}{2} + O(D), r - O(\sqrt{D})
\right)$, as $D\to 0$. This is confirmed by the numerical results, as can be seen in the final panel of Figure~\ref{fig2}. This completes the asymptotic structure for steady states $D$ small and $a>\frac{1}{2} + o(1).$\\\\
We now move on to the final case, with $D=O(1)$ and $a>\pi\sqrt{D}$. 

\subsubsection{$D=O(1)$ with $a>\pi\sqrt{D}$}
\begin{figure}
\begin{center}
\includegraphics[width=0.75\textwidth]{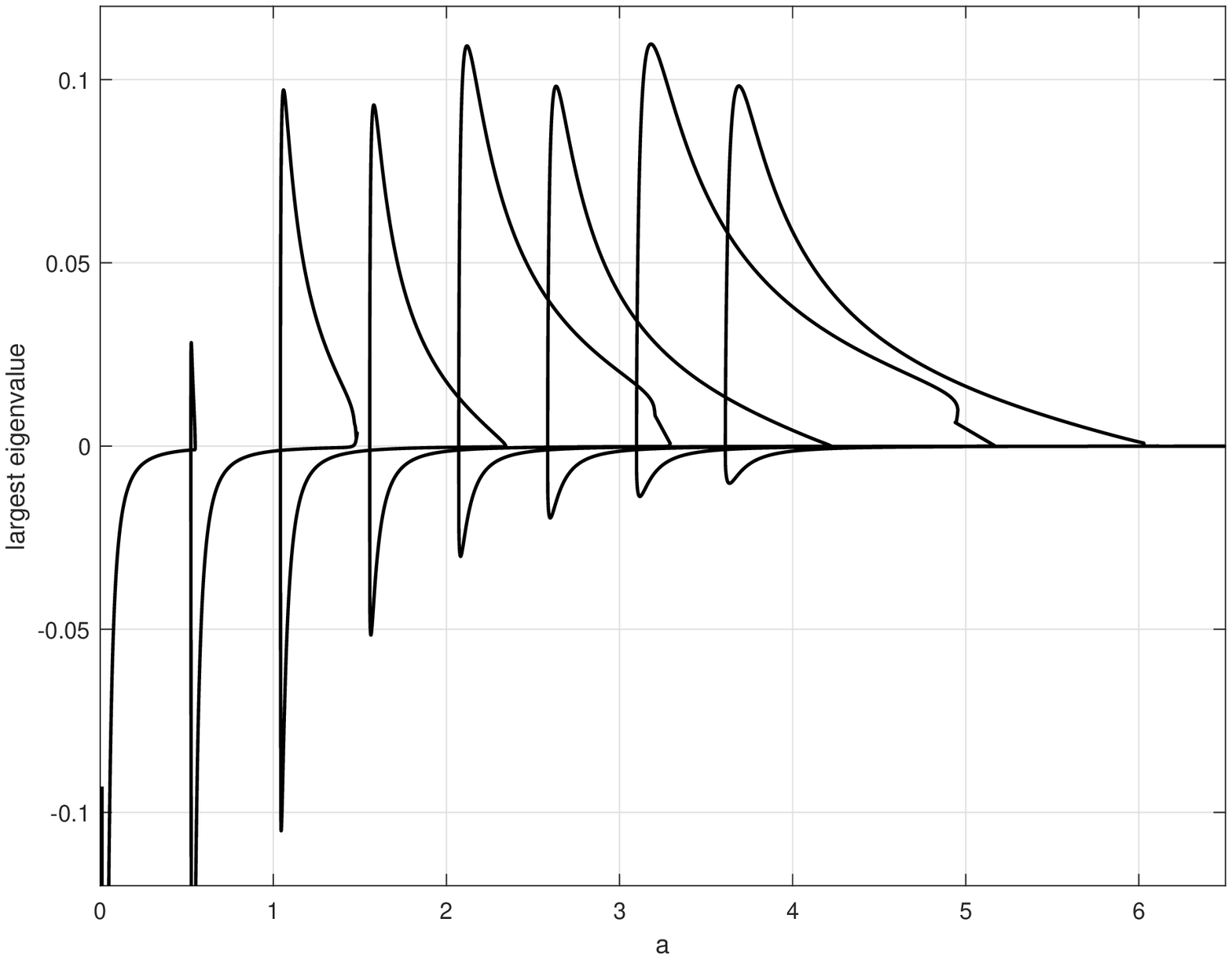}
\caption{The largest eigenvalue of the linear stability problem for $D = 10^{-5}$ as a function of $a$ as the bifurcation curve is traversed. Note the slight qualitative difference in the plot between steady states with odd and even numbers of peaks.}\label{D1em5_ev}
\end{center}
\end{figure}

\begin{figure}
\begin{center}
\includegraphics[width=0.75\textwidth]{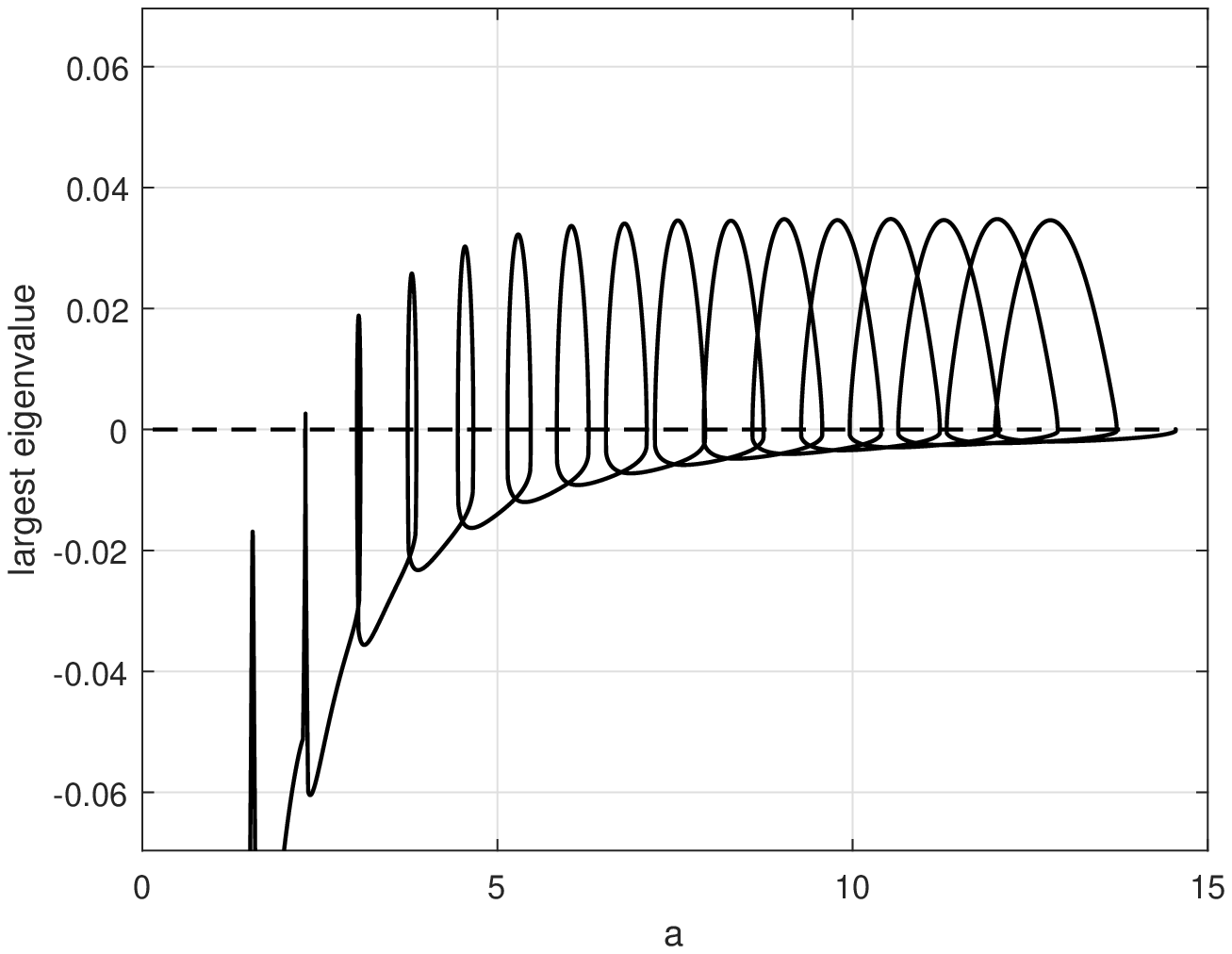}
\caption{The largest eigenvalue of the linear stability problem for $D = 2 \times 10^{-3}$ as a function of $a$ as the bifurcation curve is traversed.}\label{D2em3_ev}
\end{center}
\end{figure}

In this final case, we first investigate numerically, using pseudo-arclength continuation. This gives strong support for us to conclude the following about the bifurcation structure in the $(a,||u_s||_1)$ bifurcation plane, for nontrivial, non-negative, steady states : 
\begin{itemize}
\item For each $D\ge \Delta_1\approx0.00297$, the branch of single peak steady states created at the transcritical bifurcation from unity, when $a=\pi\sqrt{D}$, continues as a monotone increasing curve for all increasing $a$, with the associated steady state retaining a single peak structure for $D\ge\Delta_m$, but developing boundary oscillations about unity which decay exponentially into the core for $\Delta_1\le D < \Delta_m$. In each case $||u_s||_1\sim a$ as $a\to \infty$, in agreement with section 2.2.4 above. There are no other non-negative steady state branches. The upper panels in Figure~\ref{fig2} show typical bifurcation curves.
\item For each $0< D < \Delta_1$, as $a$ increases in the $(a,||u_s||_1)$ bifurcation plane, first there is the single peak steady state which emerges from the transcritical bifurcation from unity, when $a=\pi\sqrt{D}$. Thereafter, for each $r=2,3,4,...$, there is an interval $K_r(D)=[a_r(D),b_r(D)]$, with $\{a_r(D)\},~\{b_r(D)\}$ both being positive, strictly increasing and unbounded sequences in $r=2,3,...$, such that for each $a\in K_r(D)$ there is an $r$-peak nontrivial and non-negative steady state, and this has $||u_s||_1 = O(r)$ on $K_r(D)$. It follows from the theory  above in section 2.2.4, for $a$ large with $D=O(1)$, that $a_r(D)\sim \lambda_1^-(D)r$ and $b_r(D)\sim \lambda_1^+(D)r$ as $r\to \infty$. In addition, there is a natural number $r_c(D) = \lfloor \lambda_1^-(D)/(\lambda_1^+(D)-\lambda_1^-(D)) \rfloor +1 \ge 2$, which is decreasing as $D$ decreases, with $r_c(0^+) = 2$ and $r_c(D)\to \infty$ as $D \to \Delta_1^-$, and for which $a_{r+1}(D)>b_r(D)$ when $2\le r < r_c(D)$ whilst $a_{r+1}(D)\le b_r(D)$ for $r\ge r_c(D)$. When $2\le r< r_c(D)$, the connection between the $r^{th}$ and $(r+1)^{th}$ echeloned branch is a monotone increasing curve in the $(a,||u_s||_1)$ bifurcation diagram, with the associated steady state smoothly transitioning from an $r$-peak to an $(r+1)$-peak steady state. However, for $r\ge r_c(D)$ the lower end of the $(r+1)^{th}$ echelon now overlaps the upper end of the $r^{th}$ echelon, and the connection is now via a loop, which has the qualitative character of the loops visible in the numerical examples illustrated in Figure~\ref{fig2}. As a particular illustration of this case Figure~\ref{D2em3_bif} shows the $(a,||u_s||_1)$ bifurcation diagram for $D=2\times 10^{-3}$, obtained numerically using pseudo-arclength continuation. For the same case, Figure~\ref{D2em3_sols} examines in detail the $r=14$ echelon, with the transitional loop connecting to the $r=13$ echelon. Associated steady state profiles are presented at points along the upper echelon and on the loop.
\end{itemize}
Each situation has now been considered. A point of interest to note, is that the appearance of the multiple saddle-node loops in the $(a,||u_s||_1)$ bifurcation diagram, at each fixed $D$ below $\Delta_1\approx0.00297$, will generate \emph{multiple hysteresis} behaviour in the evolution of stable steady states when the parameter $a$ is slowly increased and then slowly decreased.\\[2\baselineskip]
This completes the study of nontrivial, non-negative steady states for $a>1/2$. We are now in a position to return to (DIVP) at parameter values $(a,D)$ with $a>1/2$ and $0<D<\pi^2/a^2$, which we label as region $\mathcal{H}$. We begin by addressing the local temporal stability properties of the nontrivial, non-negative steady states at each point $(a,D)\in \mathcal{H}$, recalling from above that at each such point there is either, 
\begin{itemize}
\item a single $r$-peak, nontrivial, non-negative steady state, for some $r\in \mathbb{N}$, when we label $(a,D)\in \mathcal{H}_1$, or,
\item an odd number, larger than unity, of  nontrivial, non-negative steady states, these being an $r$-peak, an $(r+1)$-peak and a transitional steady state, from a minimum value of $r(\ge 2)$ to a maximum value of $r$ ( which becomes unbounded with increasing $a$ at fixed $D<\Delta_1$ )  when we label $(a,D)\in \mathcal{H}_2$. For $D$ small the minimum and maximum values of $r$ are well approximated by the integer parts of $a+\frac{1}{2}$ and $2a+1$ respectively.
\end{itemize}
It is again helpful to fix $D$ and then follow the nontrivial and non-negative steady state bifurcation curve in the $(a,||u_s||_1)$ bifurcation plane. As this curve is followed from the initial transcritical bifurcation at $a=\pi\sqrt{D}$ to large $a$, computational results determine that the largest eigenvalue of the linear stability problem (see Appendix~\ref{app_stab}), which is zero at the small norm transcritical bifurcation, becomes negative with increasing $a$, and remains negative, unless a fold is traversed, in which case it changes sign to positive, returning to negative after the next fold is traversed, with this pattern continuing as the bifurcation curve is continuously traversed. This enables us to conclude that at any point $(a,D)\in \mathcal{H}_1$, then the $r$-peak steady state is linearly temporally stable, whilst when $(a,D)\in \mathcal{H}_3$, for each r, both the $r$-peak and the $(r+1)$-peak steady states are linearly temporally stable, with the transitional steady state being linearly temporally unstable. Figures~\ref{D1em5_ev} and~\ref{D2em3_ev} show the value of the largest eigenvalue as a function of $a$, at each point on the bifurcation curve, at two representative values of $D$, namely $D=10^{-5}$ and $D=2\times10^{-3}$, respectively. 
The eigenvalues are determined by pseudo-arclength continuation following the steady state bifurcation curve, and formulating and numerically solving the associated self-adjoint linear eigenvalue problem. The points which have zero largest eigenvalue correspond to the upper and lower saddle-node bifurcations on each echelon-loop-echelon cycle, with the largest eigenvalue being negative on the echelons, and positive on the loops. 

\begin{sidewaysfigure}
\begin{center}
\vspace{12.5cm}
\includegraphics[width=\columnwidth]{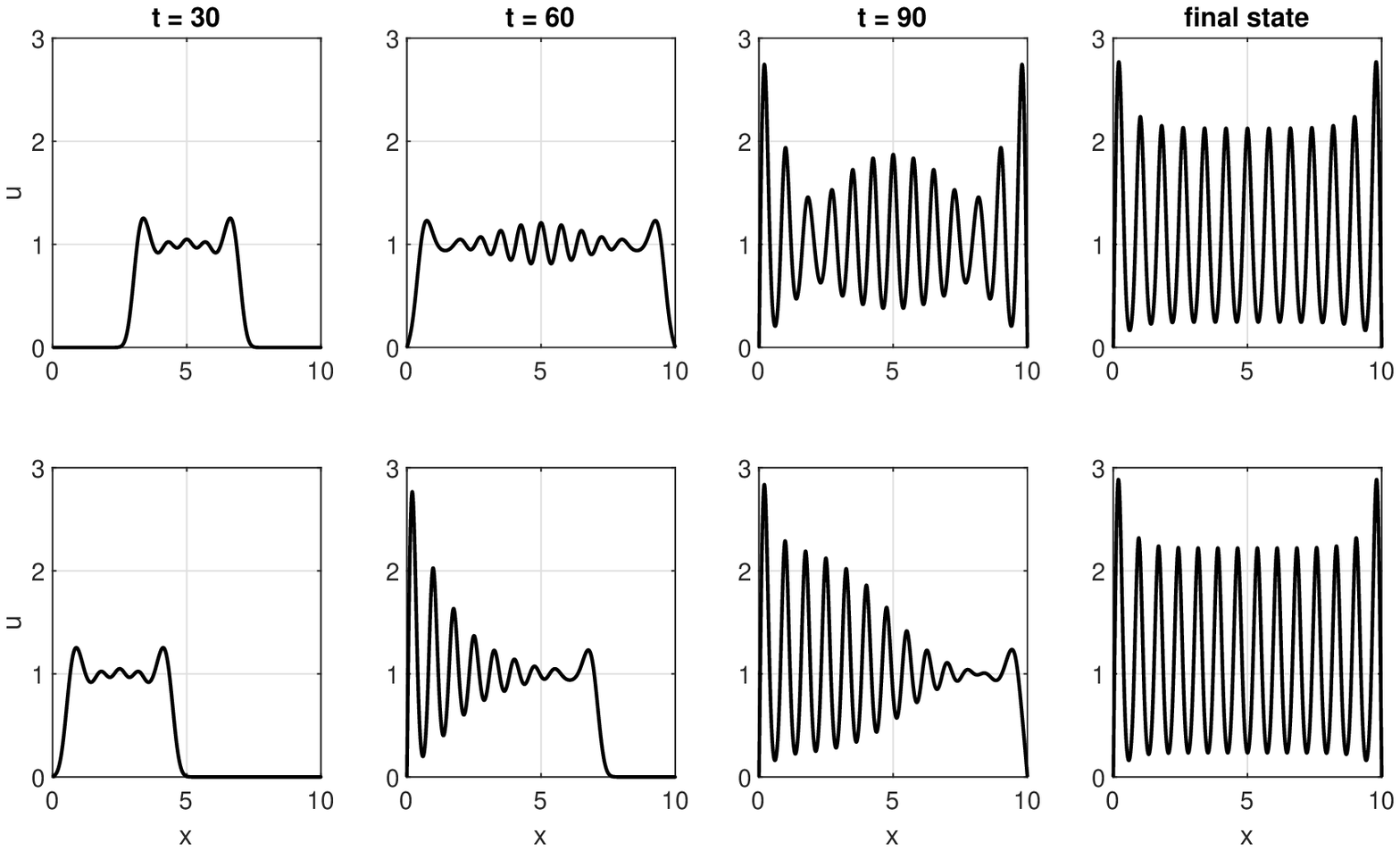}
\caption{Solutions of (DIVP) for two different initial conditions, with $a=10$ and $D = 2 \times 10^{-3}$. In the top row, the initial small input of $u$ lies symmetrically at $x_0 = 5$ and generates a steady state with 13 peaks. In the bottom row, the initial small input of $u$ lies at $x_0 = 2.5$ and generates a steady state with 14 peaks. In each case $w = 0.1$ and $\alpha = 0.01$.}\label{fig_IVP}
\end{center}
\end{sidewaysfigure}

Interpreting this in relation to (DIVP), we conjecture that :
\begin{itemize}
\item For $(a,D)\in \mathcal{H}_1$, the solution to (DIVP) approaches the temporally stable unique $r$-peak steady state as $t\to \infty$ uniformly for $x\in [0,a]$.
\item For $(a,D)\in \mathcal{H}_2$, the solution to (DIVP) approaches one of either the temporally stable $r$-peak or the temporally stable $(r+1)$-peak steady state, for some $r$ between the available minimum and maximum, depending upon the initial data $u_0:[0,a]\to \mathbb{R}$, as $t\to \infty$ uniformly for $x\in [0,a]
$.
\end{itemize}
To support this conjecture, we have solved (DIVP) numerically (see Appendix~\ref{app_num} for details) for initial conditions of the form
\begin{equation}
    u(x,0) = \left\{
    \begin{array}{cc}
    \alpha \left(1-\frac{4\left(x-x_0\right)^2}{w^2}\right) & \mbox{for $|x-x_0|<\frac{1}{2}w$},\\
    0 & \mbox{for $|x-x_0| \geq \frac{1}{2}w$.}
    \end{array}\right.
\end{equation}
This is a localised initial input of $u$ based at $x_0$, with width $w$ and height $\alpha$, each chosen so that the support of $u(x,0)$ is contained in $[0,a]$. For $D=2\times10^{-3}$ and $a=10$, two numerical solutions, with small, localised inputs of $u$ ($\alpha = 0.01$, $w = 0.1$) are shown in Figure~\ref{fig_IVP}. 
In each case, initially a pair of wavefronts is generated, which propagate away from $x = x_0$. In the symmetric case (top row), $x_0 = 5$, these reach the boundaries simultaneously and a steady state with $13$ peaks is generated. In the bottom row,  $x_0 = 2.5$, and the wavefronts reach the boundaries at different times. In this case, a steady state with $14$ peaks is generated. Figure~\ref{D2em3_bif} shows that when $a=10$, there are actually three temporally stable steady states available, with $r = 13$, $14$ and $15$. Whatever value we choose for $x_0$, numerical solutions indicate that the $r=15$ steady state is not generated using these initial conditions, and we hypothesise that this is because the slightly shorter wavelength of the $r=15$ steady state is too far from the natural wavelength generated behind the wavefronts when the Cauchy problem is on the whole real line, as determined in (NB1).

The analysis of the Dirichlet problem (DIVP) is now complete, and we move on to consider the Neumann problem (NIVP).

\section{The Neumann problem (NIVP)}
In this section we analyse the Neumann problem (NIVP). The analysis is very similar to that in section 2 for the Dirichlet problem (DIVP), and so we focus mainly on presenting the key results, avoiding repetition where possible.
When the initial data is trivial the solution to the associated (NIVP) is again the equilibrium solution given in (\ref{eqn2.1}). For initial data with $||u_0||_{\infty}$ small, it is again straightforward to develop a linearised theory for (NIVP), and this establishes that the trivial equilibrium solution is temporally unstable for all $(a,D)\in (0,\infty)\times(0,\infty)$.
\subsection{$0<a\le 1/2$}

We now consider (NIVP) when $0<a\le 1/2$, and so equation (\ref{eqn1.1}) takes on the form of (\ref{eqn2.2}) with (\ref{eqn2.3}). Firstly, we can readily establish that the only nontrivial steady state which exists is the constant steady state,
\begin{equation}
    u_s(x,D,a) = a^{-1}~~\forall~~(x,t)\in \bar{D}_{\infty}. \label{eqn3.1}
\end{equation}
Moreover, the solution to (NIVP) can be determined exactly, and is given by,
\begin{equation}
    u(x,t) = N(t)\left(\frac{1}{2}c_0 + \sum_{n=1}^{\infty}c_n e^{-\frac{n^2\pi^2D}{a^2}t}\cos\left(\frac{n\pi x}{a}\right)\right)~~\forall~~(x,t)\in \bar{D}_{\infty}  \label{eqn3.2}
\end{equation}
with
\begin{equation}
    N(t) = \left( \frac{1}{2}ac_0 + (1-\frac{1}{2}ac_0)e^{-t}\right)^{-1}~~\forall~~t\in (0,\infty) \label{eqn3.3}
\end{equation}
and
\begin{equation}
    c_r = \frac{2}{a}\int_0^a{u_0(y)\cos\left(\frac{r\pi y}{a}\right)}dy~~\text{for}~~r=0,1,2... .  \label{eqn3.4}
\end{equation}
We observe that for all initial data,
\begin{equation}
    u(x,t) \to a^{-1}~~\text{as}~~t\to \infty,  \label{eqn3.5}
\end{equation}
uniformly for $x\in [0,a]$. This establishes that the constant steady state $a^{-1}$ is globally temporally asymptotically stable. This case is now complete.

\subsection{$a>1/2$}
Here we consider (NIVP) when $a>1/2$. We begin by considering the possible nontrivial and non-negative steady states of equation (\ref{eqn1.1}). Again, it is convenient to fix $D>0$ and increase $a>1/2$, starting by continuing from the unique constant steady state $u_s = 2$, when $a=1/2$. There are a number of cases where we can develop an asymptotic theory for steady states, and we deal with these first.  We begin by considering this when $D$ is large.

\subsubsection{$D\gg1$}
There are two cases. First we address the case when $a\in (1/2,o(\sqrt{D}))$ as $D\to \infty$, and expand in the form
\begin{equation}
u_s(x,D,a) = \bar{u}_0(x,a) + D^{-1}\bar{u}_1(x,a) + O(D^{-2})~~\text{as}~~D\to \infty,  \label{eqn3.6}
\end{equation}
uniformly for $x\in [0,a]$. After substitution into equation (\ref{eqn1.1}) and boundary conditions (\ref{eqn1.7}) we find that 
\begin{equation}
\bar{u}_0(x,a) = c(a)~~\forall~~x\in[0,a],   \label{eqn3.7}
\end{equation}
with $c(a)$ a positive constant to be determined. At next order we obtain the following linear, inhomogeneous boundary value problem for $\bar{u}_1$, namely,
\begin{equation}
\bar{u}''_{1} = c(a)\left(1 - c(a)\phi(x,a)\right),~~x\in (0,a),
\end{equation}
\begin{equation}
\bar{u}'_{1}(0,a)=\bar{u}'_1(a,a)=0.
\end{equation}
Here, for $a\in (1/2,1)$,
\begin{equation}
\phi(x,a)=
\begin{cases}
x+\frac{1}{2},~~x\in [0,a-\frac{1}{2}],\\
a,~~x\in (a-\frac{1}{2},\frac{1}{2}), \\
a+\frac{1}{2}-x,~~x\in [\frac{1}{2}, a],
\end{cases}
\label{eqn3.8}
\end{equation}
whilst for $a\in [1,\infty)$,
\begin{equation}
\phi(x,a)=
\begin{cases}
x+\frac{1}{2},~~x\in [0,\frac{1}{2}],\\
1,~~x\in (\frac{1}{2},a-\frac{1}{2}), \\
a+\frac{1}{2}-x,~~x\in [a-\frac{1}{2}, a].
\end{cases}
\label{eqn3.9}
\end{equation}
The solution to this boundary value problem has,
\begin{equation}
\bar{u}'_1(x,a) = c(a)\left(x - c(a)\int_0^x{\phi(s,a)}ds\right)~~\forall~~x\in [0,a],
\end{equation}
with now $c(a)$ being required to satisfy the solvability condition,
\begin{equation}
c(a)\int_0^a{\phi(s,a)}ds - a = 0. \label{eqn3.11}
\end{equation}
Equation (\ref{eqn3.11}) has the unique positive solution,
\begin{equation}
c(a) = \frac{a}{(a-\frac{1}{4})}.  \label{eqn3.12}
\end{equation}
Thus for each $a\in (1/2,o(\sqrt{D}))$ there is a unique, nontrivial and non-negative steady state, which has,
\begin{equation}
u_s(x,D,a) = \frac{a}{(a-\frac{1}{4})}~ +~D^{-1}\left(c(a)\int_0^x{\left(y - c(a)\int_0^y{\phi(s,a)}ds\right)}dy + d(a)\right)~ +~O(D^{-2}),  \label{eqn3.13}
\end{equation}
as $D\to \infty$, uniformly for $x\in [0,a]$, with the constant $d(a)$ determined via a solvability condition at $O(D^{-2})$. It should be noted that this steady state is, as it should be, an even function about $x=\frac{1}{2}a$.  We now observe that expansion (\ref{eqn3.13}) becomes nonuniform when $a = O(\sqrt{D})$ as $D\to \infty$. Therefore, to continue this unique branch of steady states for $a$ large, we first write,
\begin{equation}
a = \sqrt{D} \bar{a}~~\text{and}~~x=\sqrt{D}\bar{x},   \label{eqn3.14}
\end{equation}
with $\bar{a}=O(1)^+$ and $\bar{x}\in [0,\bar{a}]$ as $D\to \infty$. It then follows from (\ref{eqn3.12}) and (\ref{eqn3.13}) that we should expand in the form,
\begin{equation}
u_s(\bar{x},D,\bar{a}) = 1 + D^{-\frac{1}{2}}U_0(\bar{x},\bar{a})  + O(D^{-1})~~\text{as}~~D\to \infty,   \label{eqn3.15}
\end{equation}
for $\bar{x}\in [0,\bar{a}]$. At leading order, the only solution, which has the required symmetry about $\bar{x}=\frac{1}{2}\bar{a}$, is given by,
\begin{equation}
U_0(\bar{x},\bar{a}) = E\cosh\left(\bar{x}-\frac{1}{2}\bar{a}\right)~~\forall~~\bar{x}\in (0,\bar{a}),  \label{eqn3.16}  
\end{equation}
with $E$ a constant to be determined. We now observe that this leading order term fails to satisfy the Neumann boundary conditions at the interval end points, and we conclude that two symmetric boundary layers are required. We will consider the boundary layer at the end $\bar{x}=0$. The boundary layer has $\bar{x} = O(D^{-\frac{1}{2}})$, and so $x=O(1)$, as $D\to \infty$, and we expand in the form,
\begin{equation}
u_s(x,D,\bar{a}) = 1 + D^{-\frac{1}{2}}{U}^0(x,\bar{a}) + D^{-1}{U}^1(x,\bar{a}) + O(D^{-\frac{3}{2}})~~\text{for}~~x\ge 0,  \label{eqn3.17}
\end{equation}
as $D\to \infty$. On substitution into equation (\ref{eqn1.1}), solving at each order, and applying the boundary condition at $x=0$, we obtain,
\begin{equation}
{U}^0(x,\bar{a}) = F,   \label{eqn3.18}
\end{equation}
\begin{equation}
{U}^1(x,\bar{a}) = 
\begin{cases}
\left(G-\frac{1}{16}-\frac{1}{8}(x-\frac{1}{2})\right),~~x>\frac{1}{2},\\
\left(G - \frac{1}{8}x - (\frac{1}{2}-x)^3\right),~~0\le x\le \frac{1}{2},
\end{cases}
\label{eqn3.19}
\end{equation}
with $F$ and $G$ constants to be determined via matching to the core region.  Matching (via Van Dyke's Principle) expansion (\ref{eqn3.15}) (as $\bar{x}\to 0$) with expansion (\ref{eqn3.17}) (as $x\to \infty$) determines,
\begin{equation}
E = \frac{1}{8~\text{sinh}\left(\frac{1}{2}\bar{a}\right)}~~\text{and}~~F= \frac{1}{8~\text{tanh}\left(\frac{1}{2}\bar{a}\right)}.  \label{eqn3.20}
\end{equation}
However, the determination of the constant $G$ requires the next term in the core region to be found, which, for brevity, we do not pursue further. 

\begin{figure}
\begin{center}
\includegraphics[width=0.8\textwidth]{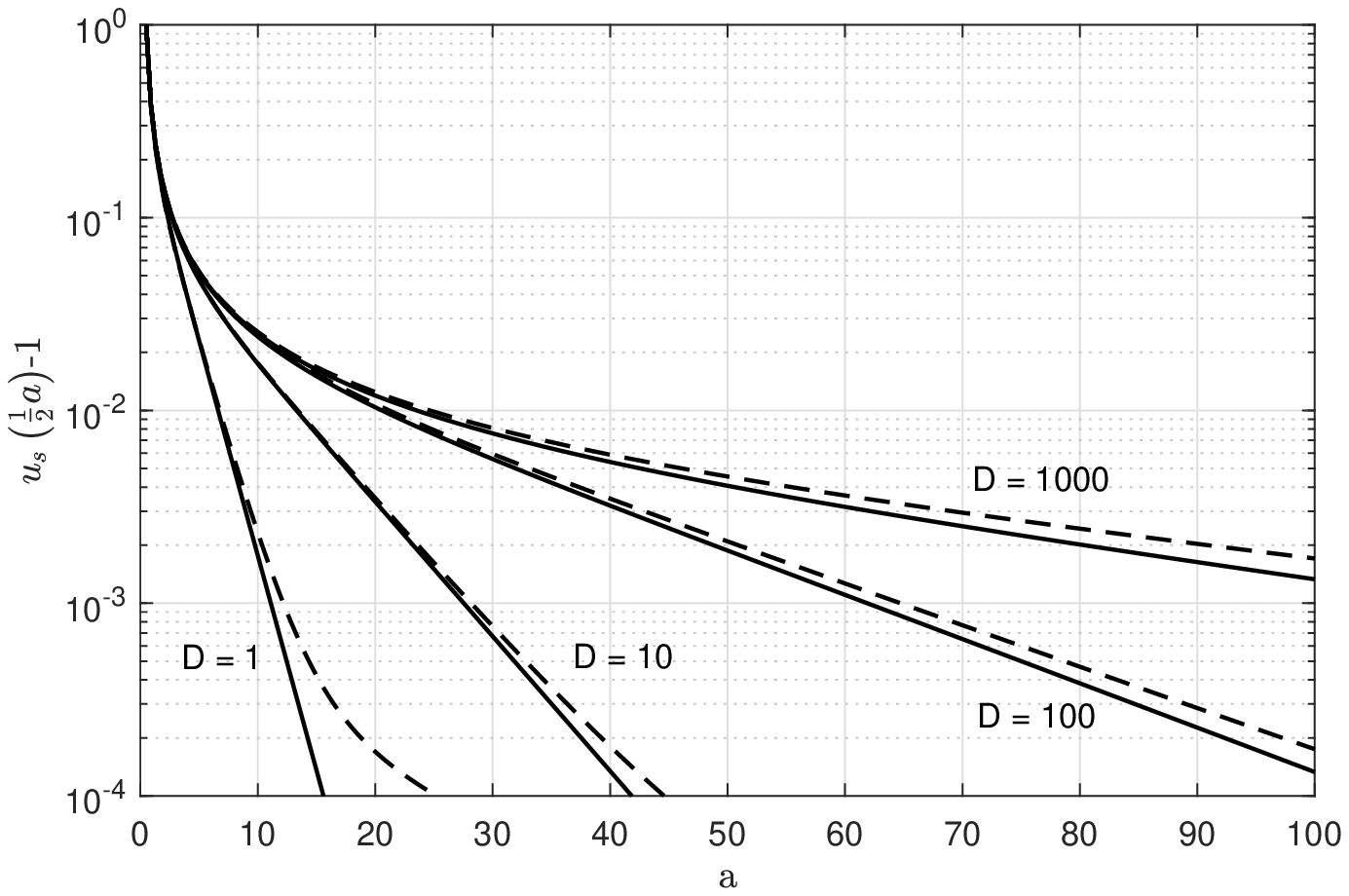}
\caption{The value of the steady state at the midpoint of the domain (minus one for clarity), $u_s\left(\frac{1}{2}a\right) - 1$, as a function of $a>\frac{1}{2}$ for $D = 1$, $10$, $100$ and $1000$. The broken line is the leading order composite asymptotic approximation for $D \gg 1$, (\ref{compos}).}\label{N_largeD}
\end{center}
\end{figure}

We now make a comparison between the above asymptotic forms for steady states valid for $D \gg 1$ and the numerical solutions obtained by pseudo-arclength continuation following steady states  of the full problem. This is done most conveniently by plotting the value of the steady state at the midpoint of the domain (minus one for clarity in a semi-logarithmic plot), $u_s\left(\frac{1}{2}a\right) -1$, as a function of $a>\frac{1}{2}$ for various values of $D$, and compare with the leading order composite form for $u_s\left(\frac{1}{2}a\right)$, given by
\begin{equation}
    u_s\left(\frac{1}{2}a,D,a\right) = \frac{a}{(a-\frac{1}{4})}+\frac{1}{8 \sqrt{D} \sinh\left(\frac{a}{2\sqrt{D}}\right)}-\frac{1}{4a} + O(D^{-1}),\label{compos}
\end{equation}
 as $D\to 0$, and obtained from the matched asymptotic expansions (\ref{eqn3.13}) and (\ref{eqn3.15}) above. This gives the leading order approximation \emph{uniformly} for all $a>\frac{1}{2}$ when $D \gg 1$. Figure~\ref{N_largeD} shows how the agreement between the numerically determined steady state and leading order composite asymptotic form of the steady state improves as $D$ increases.

This completes the analysis of the steady state structure when $D$ is large. We now consider the situation when $D=O(1)$ and $a$ is large.

\subsubsection{$D=O(1)$ with $a\gg1$}
In this case the details are very similar to the corresponding case in section 2, as detailed in subsection 2.2.4, and so we limit ourselves to presenting the salient results. For $D\ge \Delta_1$, there is a unique steady state at each large $a$, which has, in the bulk region,
\begin{equation}
u_s(x,D,a)\sim 1~~\text{as}~~a\to \infty
\end{equation}
for $x\in (O(1)^+,a-O(1)^+)$. However, for $0<D< \Delta_1$ the existence of multi-peak steady states appears for $a$ large. In particular, for fixed $D<\Delta_1$ , and each large integer $r$, there is an $r$-peak nontrivial, non-negative steady state for each $\lambda\in (\lambda_1^-(D),\lambda_1^+(D))$, given by, in the bulk region,
\begin{equation}
u_s(x,D,a)\sim F_p(x+x_M(\lambda,D),\lambda,D),  \label{eqn3.21}
\end{equation}
for $x\in (O(1)^+,a-O(1)^+)$, with corresponding $L^1$ norm being,
\begin{equation}
||u_s(\cdot,D,a)||_1\sim rA(\lambda,D)  \label{eqn3.22} 
\end{equation}
 and large $a$ given, in terms of $\lambda$, by
\begin{equation}
a = \lambda(r-1) + O(1)~~\text{as}~~r\to \infty,  \label{eqn3.23}
\end{equation}
for each $\lambda\in (\lambda_1^-(D),\lambda_1^+(D))$. In the $(a,||u_s(\cdot,D,a)||_1)$ bifurcation plane, for each large integer $r$, the associated family of $r$-peak steady states lies on an echelon above the interval with $a\in I_r(D)$, parameterised by $\lambda\in (\lambda_1^-(D),\lambda_1^+(D))$, via (\ref{eqn3.22}) and (\ref{eqn3.23}). Sequential echelons are again overlapping, and at a given large $a$, the overlapping echelons have those integer values, with $r_m(a,D)\le r \le r_M(a,D)$, as for the Dirichlet case. Also, for $r$ large, each echelon has length of $O(r)$ and variation of $O(r)$, and so slope of O(1), whilst each sequential echelon is lifted by $O(1)$ with each unit increase in $r$. Again, in the bifurcation diagram, each sequential echelon is connected by a loop containing a supercritical and subcritical saddle-node bifurcation point, the traversal of which accommodates the transition from an $r$-peak to an $(r+1)$-peak structure in the steady state. We next consider the case when $D$ is small. 

\subsubsection{$0<D\ll1$}
\begin{figure}
\begin{center}
\includegraphics[width=0.8\textwidth]{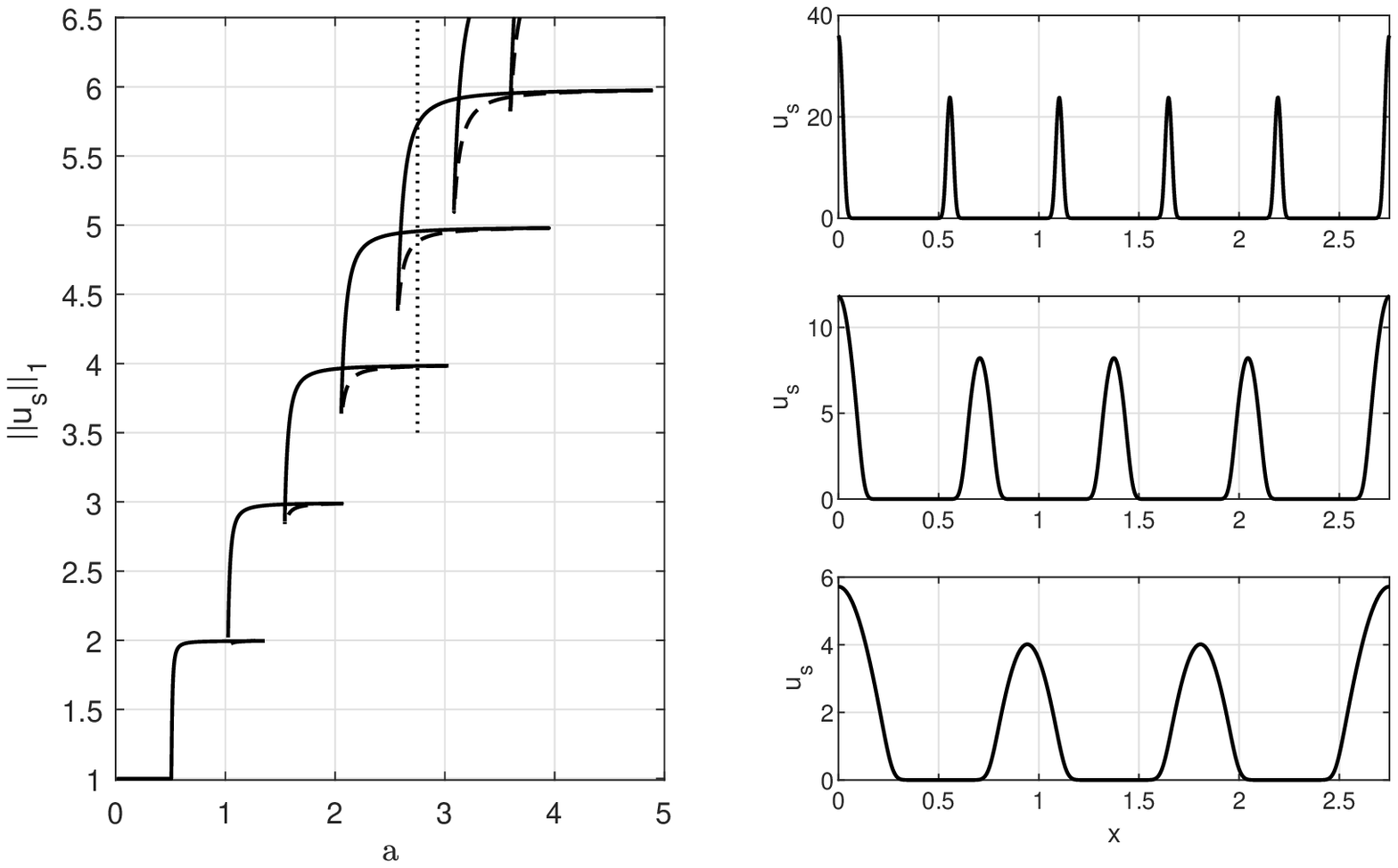}
\caption{The bifurcation diagram for (NIVP) with $D = 10^{-5}$ (left panel). The three stable steady states when $a = 2.75$ (indicated by the dotted line), with $r = 4$, $5$ and $6$, are shown in the right hand panel to illustrate the typical form of these steady states. Note the different scales on the $u_s$-axes.}\label{fig_UNeu}
\end{center}
\end{figure}

The case $0<D\ll1$ is similar to the corresponding case for the Dirichlet problem dealt with in section 2, and we therefore present only the salient features here. For each $r=2,3,...$ we can construct a branch of $r$-peak nontrivial, non-negative steady states, for
\begin{equation}
    a\in \left(\frac{1}{2}(r-1),r-\frac{1}{2}\left(3-\sqrt{2}\right)\right),    \label{eqn3.24}
\end{equation}
given by,
\begin{equation}
    u_s(x,D,a) \sim 
    \begin{cases}
   \frac{\pi}{\sqrt{2}(\lambda(a,r)-\frac{1}{2})}\cos\left(\frac{\pi x}{\sqrt{2}(\lambda(a,r)-\frac{1}{2})}\right),~~0\le x \le \frac{1}{\sqrt{2}}\left(\lambda(a,r)-\frac{1}{2}\right),\\
    F_p\left(x-\frac{1}{\sqrt{2}}(\lambda(a,r)-\frac{1}{2}),\lambda(a,r),D\right),~~\frac{1}{\sqrt{2}}\left(\lambda(a,r)-\frac{1}{2}\right)<x<a-\frac{1}{\sqrt{2}}\left(\lambda(a,r)-\frac{1}{2}\right),\\
   \frac{\pi}{\sqrt{2}(\lambda(a,r)-\frac{1}{2})}\cos\left(\frac{\pi (a-x)}{\sqrt{2}(\lambda(a,r)-\frac{1}{2})}\right) ,~~a-\frac{1}{\sqrt{2}}\left(\lambda(a,r)-\frac{1}{2}\right)\le x \le a,
    \end{cases}
    \label{eqn3.25}
\end{equation}
as $D\to 0$, with
\begin{equation}
  \lambda(a,r) = \left(a + \frac{1}{2}\left(\sqrt{2}-1\right)\right) (r - 2 + \sqrt{2})^{-1}.
\end{equation}
Here, following (NB1), we have over one wavelength $\lambda$, with $\lambda\in \left(\frac{1}{2},1\right)$,
\begin{equation}
    F_p(y,\lambda,D) \sim 
 \begin{cases}

0,~~~0 \le y < \frac{1}{2}, \\
\frac{\pi}{2(\lambda-\frac{1}{2})}\cos\left(\frac{\pi(y-\frac{1}{2}(\lambda+\frac{1}{2}))}{(\lambda-\frac{1}{2})}\right),~~~ \frac{1}{2}\le y\le \lambda, 
\end{cases}  \label{eqn3.26}
\end{equation}
as $D\to 0$. Again, localised edge layers, of thickness $O(D^{\frac{1}{4}})$, are present to smooth $F_p$ across the edge of its support regions in (\ref{eqn3.26}) (see section 5.1  in (NB1)). A notable feature of the $r$-peak steady state in the Neumann case, given by (\ref{eqn3.25}), is that the peaks occuring at both end points have height and support width which is exactly $\sqrt{2}$ times that of each interior peak. This is due to one half of the spatial extent of the boundary peak being clipped off beyond the boundary in relation to the nonlocal term in equation (\ref{eqn1.1}), whilst the edge layer requires that the slope at the edge of the boundary peak is equal to that at the edges of the interior peaks. This does not occur in the Dirichlet case, when all peaks have the same height in this limit. In the $(a,||u_s||_1)$ bifurcation plane the $r$-peak echelon, for $r=2$, $3$, $\ldots$, has the curve given by,
\begin{equation}
 ||u_s(\cdot,D,a)||_1 = r\left(1 - \frac{\pi^2}{\left(\lambda(a,r)-\frac{1}{2}\right)^2}D + o(D^{\frac{5}{4}})\right)~~\text{as}~~D\to 0   \label{eqn3.27}
\end{equation}
for $a\in \left(\frac{1}{2}(r-1),(r-\frac{1}{2}\left(3-\sqrt{2}\right))\right)$. Finally, at each $a>1/2$ there is an $r$-peak steady state for each
\begin{equation}
r \in [1+a,1+2a].   \label{eqn3.28}
\end{equation}
As in the Dirichlet case, on the $(a,||u_s||_1)$ bifurcation plane, each sequential echelon is connected by an elongated loop which has a low curvature supercritical saddle-node bifurcation at the lower end, and a high curvature subcritical saddle-node bifurcation at the upper end. A numerical illustration of the bifurcation plane, at $D=10^{-5}$, is given in Figure~\ref{fig_UNeu}, along with the three steady state solutions when $a = 2.75$. 

This completes the analysis of nontrivial, non-negative steady states for small $D$.

\subsubsection{$D=O(1)$ with $a>1/2$}
Steady states in this final case are analysed via numerical pseudo-arclength continuation by fixing $D$ and then path following from the unique nontrivial steady state at $a=1/2$ with increasing $a$. The results are not qualitatively different to those that we obtained for the Dirichlet problem, and we have not added any new illustrations. The key  conclusions from this numerical investigation follow directly those of subsection 2.2.6 for the Dirichlet case, and are therefore not repeated again in detail. \\ \\
This completes the analysis of nontrivial, non-negative steady states for the Neumann problem when $a>1/2$. We note that the remarks concerning \emph{steady state hysteresis} in section 2 relating to the Dirichlet case, apply equally here. We are now in a position to consider the evolution problem (NIVP). At each $(a,D)\in (1/2,\infty)\times(0,\infty)$, numerical solutions to (NIVP) approach one of the temporally stable nontrivial, non-negative steady states when $t$ is large. As we saw for (DIVP), different initial conditions lead to different final steady states, but the behaviour of unsteady solutions to (NIVP) are not qualitatively different to those of (DIVP), and we do not present any further numerical illustrations.

\section{Conclusions}
In this paper, we considered two evolution problems, labeled as (DIVP), with Dirichlet boundary conditions, and (NIVP), with Neumann boundary conditions. The key feature of the inclusion of the nonlocal term is the introduction of a nonlocal length scale into the model. There are thus two dimensionless parameters in the model, namely $D$, which measures the square of the ratio of the diffusion length scale (based on the kinetic time scale) to the nonlocal length scale, and $a$, which measures the ratio of the domain length scale to the nonlocal length scale. We have examined both the Dirichlet and the Neumann models for $(a,D)$ throughout the positive quadrant of the parameter plane. Significant structural differences in behaviour between the classical Fisher-KPP model and this natural nonlocal extension have been identified, and these become more significant as the parameter $D$ decreases. These principal, and marked, differences have been fully explored and exposed. The structural differences have included the existence of multiple, temporally stable, multipeak, nontrivial and non-negative steady states, accumulated through multiple steady state saddle-node bifurcations, with this giving rise to hysteretic behaviour with increasing and decreasing domain size $a$, at sufficiently small $D$. In addition, the temporally stable multipeak steady states, when available, are shown to be large-$t$ attractors for the evolution problems (DIVP) and (NIVP), respectively.

 More generally, this paper and (NB1) have both studied a natural extension of the classical Fisher-KPP model, with the introduction of the simplest possible nonlocal effect into the saturation term. The motivation for these detailed analyses arises since nonlocal reaction-diffusion models occur naturally in a variety of (frequently biological or ecological) contexts, and as such it is of fundamental interest to examine the properties of the nonlocal Fisher-KPP model in detail, and to compare and contrast these with the well known properties of the classical Fisher-KPP model. As we have seen here, in (NB1) and in the work of many other authors (for example, \cite{BNPR,JBNL,Gourley2000,NADIN2011553,FAYE20152257,perthame2007}) the introduction of a nonlocal kernel leads to a huge increase in the complexity of possible solutions and their bifurcation structure. In particular, the detailed effect of choosing a kernel with compact support, as opposed to a kernel, for example a Gaussian, supported on the whole real line is not clear. Since a top hat kernel is the simplest possible compactly-supported kernel, the detailed analysis that we present here and in (NB1) forms the basis for an investigation into the form of solutions when the kernel is not constant on its support. Our ongoing investigation suggests that adding a small, variation to the kernel is, in the limit $D \to 0$, a singular perturbation, and that unravelling the bifurcation structure, starting from the structure presented here for the top hat kernel, will reveal a wide variety of solutions. This will be the subject of a subsequent paper.

\bibliographystyle{amsplain}
\bibliography{main}


\begin{appendix}
\section{Numerical methods}\label{app_num}
\subsection{Finite difference discretisation}\label{app_fd}

All of the numerical methods that we use in this paper are based on central finite differences in space for second derivatives and the trapezium rule for the convolution integral. We discretise at $N$ equally-spaced points, $x_i$, in $[0,a]$, with $x_0=0$ and $x_N = a$. The trapezium rule that we use for the convolution integral respects the fact that the end points, $\alpha(x)$ and $\beta(x)$, defined in (\ref{eqn1.3}) and (\ref{eqn1.4}), may fall between two grid points, and is calculated based on a linear variation of $u$ between grid points. This is used in the discretisation of both the steady and unsteady versions of (\ref{eqn1.1}) and the stability problem described in Appendix~\ref{app_stab}. The Neumann boundary conditions, (\ref{eqn1.7}), are discretised using a three point difference formula in order to maintain second order accuracy in space. We typically used $N=1000$.

\subsection{Pseudo-arclength continuation}
We investigate the bifurcation diagram for the steady states, $u_s(x)$, of (\ref{eqn1.1}) using pseudo-arclength continuation, \cite{arclength}, in $(a, ||u_1||_1)$ parameter space. We will use the notation $A \equiv ||u_s||_1$ below for clarity. The discretised linear system that arises from (\ref{eqn1.1}) with $u_t=0$ and either Dirichlet or Neumann boundary conditions is supplemented by the conditions
\begin{equation}
    \left(a-a_0\right)^2 + \left(A-A_0\right)^2 = ds^2,
\end{equation} 
and
\begin{equation}
    \int_0^a u_s(x) dx = A,
\end{equation}
with the trapezium rule used to evaluate this integral. Here, $a_0$ and $A_0$ are the values of the length of the domain and the area under the solution at the previous point on the solution curve, and $ds$ is the pseudo-arclength that we attempt to traverse in the $(a,A) \equiv (a, ||u_s||_1)$ parameter space. The nonlinear system of algebraic equations that determine the discretised solution is solved using 'fsolve' in Matlab, with the Jacobian of the system supplied analytically. The initial discrete solution is determined for $a < \frac{1}{2}$, in which case a simple analytical solution is available for both (DIVP) and (NIVP). Linear extrapolation is used to obtain new guesses of $a$ and $A$, and we use a simple scaling factor for $ds$ so that it is reduced when the solution fails to converge and increased otherwise, with $10^{-6} \leq ds \leq 5 \times 10^{-3}$. For $D$ sufficiently large, specifically $D>10^{-3}$, this method is sufficient to calculate the bifurcation diagrams presented in the paper, including the loops in parameter space shown in Figure~\ref{fig2}.

For smaller values of $D$, the curvature of the saddle node bifurcation at the right of each loop becomes so large that this algorithm fails. It seems likely that this curvature becomes exponentially large as $D \to 0$, in effect forming a cusp at which the direction of the bifurcation curve reverses at the bifurcation, which is too much for our simple algorithm to cope with. We overcome this by using parameter continuation to generate a solution on the stable branch above the cusp and then use pseudo-arclength continuation to reach the cusp from the opposite direction. This is sufficient to generate all of the bifurcation diagrams in the paper. At all values of $D$ for which bifurcation loops occur, when $a$ becomes sufficiently large the curvature of the saddle-node bifurcation at the left of the loop also becomes too tight for our algorithm to cope with, but we have not attempted to deal with this problem as the results presented in the paper are sufficient to illustrate the form of the bifurcation diagram.

\subsection{Numerical solution of (DIVP) and (NIVP)}
After discretising (DIVP) and (NIVP) in space using the finite difference method described in Appendix~\ref{app_fd}, we use the midpoint method to evolve the solution from its initial state. If we write the discretised system as $\dot{\bf u} = {\bf F}({\bf u})$, with $u_i = u(x_i,t)$, this is
\[
    \bar{\bf u}(t) = {\bf u}(t) + \frac{1}{2}\Delta t \, {\bf F}({\bf u}(t)),~~~{\bf u}(t+\Delta t) = {\bf u}(t) + \Delta t \, {\bf F}(\bar{\bf u}(t)).
\]
We choose a time step, $\Delta t$, adaptively so that the maximum change in $u$ is less than $10^{-4}$, with maximum time step $\frac{1}{4} D^{-1} \Delta x^2$.

\section{Linear stability of steady states}\label{app_stab}
To determine the linear stability property of the steady state $u_s$ we define a new dependent variable  $\hat{u}(x)$ by
\[u= u_s(x) + \hat{u}(x) e^{\sigma t},\]
substitute into (DIVP) and (NIVP), and assume that $|\hat{u}| \ll 1$, to obtain the linear stability problem
\begin{equation}
    D \hat{u}'' + \left(1 - \int_{\alpha(x)}^{\beta(x)} u_s(y)dy\right) \hat{u} -  u_s(x) \int_{\alpha(x)}^{\beta(x)} \hat{u}(y)dy = \sigma \hat{u},~~x\in (0,a),
\end{equation}
subject to the boundary conditions
\begin{equation}
\hat{u}(0) = \hat{u}(a) = 0,
\end{equation}
for (DIVP), and
\begin{equation}
\hat{u}'(0) = \hat{u}'(a) = 0,
\end{equation}
for (NIVP). This linear eigenvalue problem is self-adjoint, and therefore has real eigenvalues, $\sigma$. It is straightforward to discretise  this problem using the finite difference approach described in Appendix~\ref{app_fd}. The resulting matrix eigenvalue problem can then be solved in Matlab using the built-in routine 'eig'. Typical results for the largest eigenvalue as a function of $a$ for fixed $D$ are shown in Figure~\ref{D1em5_ev} and~\ref{D2em3_ev}. 
\end{appendix}

\end{document}